\documentclass[preprint]{elsarticle}
\usepackage{amsfonts,amsmath,amssymb}
\usepackage{graphicx}
\usepackage{algpseudocode}
\usepackage{algorithm}
\usepackage{multirow}

\usepackage{colortbl}
\usepackage{natbib}

\newcommand{\cond}{\,|\,}

\DeclareMathOperator{\trace}{tr}

\DeclareMathOperator{\E}{\mathbb{E}}
\DeclareMathOperator{\KL}{KL}

\begin{document}
\title{Gaussian filtering and variational approximations for Bayesian smoothing in continuous-discrete stochastic dynamic systems}
\author[tut1]{Juha Ala-Luhtala}
\author[aalto]{Simo S\"arkk\"a}
\author[tut1]{Robert Pich\'e}

\address[tut1]{Tampere University of Technology, Finland}
\address[aalto]{Aalto University, Finland}

\begin{abstract}
The Bayesian smoothing equations are generally intractable for systems described by nonlinear stochastic differential equations and discrete-time measurements. 
Gaussian approximations are a computationally efficient way to approximate the true smoothing distribution. 
In this work, we present a comparison between two Gaussian approximation methods.
The Gaussian filtering based Gaussian smoother uses a Gaussian approximation for the filtering distribution to form an approximation for the smoothing distribution. 
The variational Gaussian smoother is based on minimizing the Kullback--Leibler divergence of the approximate smoothing distribution with respect to the true distribution. 
The results suggest that for highly nonlinear systems, the variational Gaussian smoother can be used to iteratively improve the Gaussian filtering based smoothing solution. 
We also present linearization and sigma-point methods to approximate the intractable Gaussian expectations in the Variational Gaussian smoothing equations. 
In addition, we extend the variational Gaussian smoother for certain class of systems with singular diffusion matrix.
\end{abstract}


\maketitle

\section{Introduction}

The continuous-discrete system refers to a system whose process dynamics are governed by a continuous-time stochastic differential equation (SDE) and whose measurements are taken at discrete time instants. 
Bayesian filtering and smoothing equations give the solution to the problem of estimating the state of the system from the noisy measurements. 
Computing the filtering and smoothing distributions involves solving the related partial differential equations \cite{Jazwinski1970, Leondes1970},  
and is only tractable for linear-Gaussian systems (and some other special cases \cite{Daum1986}).
In this paper, we consider two different approaches for computing a Gaussian approximation for the smoothing distribution. 

The Gaussian approximation for the filtering and smoothing equations is well known in the literature. 
The first Gaussian approximations were based on linearization using the Taylor series based methods \cite{BrysonFrazier1963, Jazwinski1970}. 
The Taylor series based methods can be seen as a special case of the more general Gaussian filtering and smoothing framework, where different filters and smoothers arise based on the numerical method for computing the Gaussian expectations \cite{WuEtAl2006, SarkkaSarmavuori2013}. 
A different approach for Gaussian smoothing was considered in \cite{Archambeau2007,Archambeau2008, ArchambeauOpper2011}, where a Gaussian approximation is sought by approximating the stochastic process giving the smoothed distribution with a linear process. 
The method is based on the fixed-form variational Bayes approximation \cite{Bishop2006} and minimizes the Kullback--Leibler divergence of the approximate distribution with respect to the true distribution. 

The variational Gaussian approximation is considered further in \cite{Shen2010, Shen2012, Vrettas2010, Vrettas2011}. 
Shen \emph{et al.} \cite{Shen2010} compared the variational approximation to a Monte Carlo Markov Chain (MCMC) solution for the one dimensional double well system and found that the variational method performed comparatively to the MCMC solution when the uncertainties in the measurements were not so large as to cause the true posterior to be multimodal.
In \cite{Shen2012} the variational Gaussian smoothing solution is used as a proposal distribution in an MCMC method to improve the efficiency of the algorithm.
The variational MCMC method was found to outperform the hybrid Monte Carlo method for sparsely observed diffusion processes. 
Vrettas \emph{et al.} \cite{Vrettas2010, Vrettas2011} considered approximating the variational Gaussian smoothing equations using a radial basis function representation for the variational parameters. 
By imposing a certain structure for the variational parameter functions, the overall number of parameters to be optimized can be reduced. 

The variational Gaussian smoothing equations derived by Archambeau \emph{et al.} \cite{Archambeau2007, Archambeau2008, ArchambeauOpper2011} require a non-singular effective diffusion matrix. 
Our first goal is to extend the variational Gaussian approximation to a certain class of singular models by considering an alternative derivation based on Girsanov's theorem. 

Another problem in the variational Gaussian algorithm is the need to compute Gaussian expectations over nonlinear functions. 
Previous works on the variational approximation \cite{Archambeau2007, Archambeau2008, Shen2010, Vrettas2010, Vrettas2011} have not presented details for computing the Gaussian expectations for general nonlinear systems.
Our second goal is to extend the variational Gaussian smoothing method for general nonlinear systems by considering numerical approximations for the Gaussian expectations.
The treatment is similar to the one in \cite{SarkkaSarmavuori2013}, where Taylor series based linearization, cubature and unscented transform based sigma-point methods and Gauss--Hermite quadrature were used to compute Gaussian expectations in the Gaussian smoothers based on the classical Gaussian filtering framework. 

In addition, we provide in this paper a comparison between the variational Gaussian approximation and the Gaussian smoothers presented in \cite{SarkkaSarmavuori2013}. 
Using a suitable change of variables the Gaussian smoothing equations can be converted to a variational form similar to the variational Gaussian smoothing equations.
The computation of the Gaussian filtering based smoothing solution is numerically stable and provides good initial conditions for the variables in the variational Gaussian smoothing algorithm. 
This might help to overcome the problems reported by Vrettas \emph{et al.}  \cite{Vrettas2010} in the initialization of the variational Gaussian algorithm for high dimensional systems. 
Also, we study if the variational Gaussian smoother can be used to iteratively improve the results from the Gaussian filtering based smoother. 

The organization of this paper is as follows. 
First we present the variational Gaussian smoother and the Gaussian smoother based on the Gaussian filtering framework for the continuous-discrete system.
Here we extend the variational Gaussian smoothing equations for a certain class of singular models. 
Next we compare theoretically the Gaussian smoothers by presenting a conversion of the Gaussian filtering based smoothing equations to the variational form by a suitable change of variables. 
This also provides the initial values for the variational Gaussian smoothing algorithm. 
The problem of numerically computing the Gaussian expectations in the variational Gaussian smoothing equations is treated in the next section.
We present Taylor series linearization and sigma-point methods to approximate the Gaussian expectations in the variational Gaussian smoothing equations. 
The paper concludes with two synthetic-data examples that are used to compare the Gaussian smoothers, and also provides comparison of the different numerical methods for the computation of the Gaussian expectations. 



%
\subsection{Problem statement}
The continuous-discrete system considered in this paper is given by
\begin{align}
dx  & = f(x, t) \, dt + L(t) d \beta(t), \\
y_k & = h_k(x(t_k)) + v_k,
\label{eq:SDE}
\end{align}
where $x(t)$ is the state, $f(x(t),t)$ is the drift term, and $\beta(t)$ is a Brownian motion stochastic process with diffusion matrix $Q(t)$. 
The effective diffusion matrix for the process is given by
\begin{equation}
\Sigma(t) = L(t)Q(t)L^T(t).
\end{equation}
The initial conditions are assumed to be normally distributed $x(t_0) \sim \text{N}(m_0, P_0)$. 
The measurement noise $\{v_k\}$ is a zero mean Gaussian white noise sequence with covariance matrix $R_k$. 
The measurement noise $\{v_k\}$, process noise $\beta(t)$ and initial conditions $x_0$ are assumed to be mutually independent. 

Let $y_1,\ldots,y_K$ be the measurements taken at discrete time instants $t_1, \ldots,t_K$. 
The solution to the Bayesian smoothing problem is the posterior distribution
\begin{equation}
p(x(t) \cond y_1,\ldots,y_K), \quad t \in [t_0, t_k].
\end{equation}

In this paper, we concentrate on Gaussian approximations for the smoothing distribution. 
That is, the smoothing distribution is approximated as
\begin{equation}
p(x(t) \cond y_1,\ldots,y_K) \approx N(x(t) \cond m(t), P(t)),
\label{eq:Gaussian_approx}
\end{equation}
where $m(t)$ is the mean function and $P(t)$ is autocovariance $P(t,t')$ at $t=t'$ for the approximating distribution. 
The smoothing problem now reduces to finding the expressions for the mean and covariance functions. 

\section{Gaussian smoothing for continuous-discrete systems}

\subsection{Variational Gaussian approximation}

The variational Gaussian approximation for the continuous-discrete smoothing problem was derived by Archambeau \emph{et al.} \cite{Archambeau2007, Archambeau2008, ArchambeauOpper2011}.
The method is based on approximating the smoothing process with a linear process
\begin{equation}
dx = [-A(t)+b(t)]dt + \sqrt{\Sigma(t)} d \beta(t),  
\label{eq:sde_linear}
\end{equation}
where $A(t)$ and $b(t)$ are parameters of the approximation and $\beta(t)$ is a Brownian stochastic process with identity diffusion matrix. 
The solution to the linear SDE (\ref{eq:sde_linear}) is a Gaussian process. 
The marginal density at each time is given by $q(x(t)) = \text{N}(x(t) \cond m(t), P(t))$, where the mean and covariance are computed from the ordinary differential equations
\begin{align}
\frac{d}{dt}m(t) & = -A(t)m(t)+b(t) \label{eq:vb_mean}\\
\frac{d}{dt}P(t) & = -A(t)P(t) - P(t)A^T(t) + \Sigma(t) \label{eq:vb_cov}. 
\end{align}
The parameters $A(t)$ and $b(t)$ are computed by minimizing the Kullback--Leibler (KL) divergence of the probability law $\mathbb{Q}_X$ of the approximating process with respect to the probability law $\mathbb{P}_{X \cond Y}$ of the true smoothing process. 
The KL-divergence is given by \cite{Archambeau2007, Archambeau2008, ArchambeauOpper2011}
\begin{equation}
\KL(\mathbb{Q}_X \, || \,  \mathbb{P}_{X \cond Y})  = \int_{t_0}^{t_K} \E_q \left[e(x(t),t) + \sum_{k=1}^K u_k(x(t))\delta(t-t_k) \right] \, dt,
\label{eq:KL_div2}
\end{equation}
where
\begin{align}
e(x(t),t) & = \frac{1}{2}[f(x(t),t)+A(t)x(t)-b(t)]^T \Sigma^{-1}(t)[f(x(t), t)+A(t)x(t)-b(t)], \label{eq:e} \\
u_k(x(t)) & = \frac{1}{2}\left[y_k - h_k(x(t))\right]^TR_k^{-1}\left[y_k - h_k(x(t))\right] \label{eq:uk}. 
\end{align}
The expectations $\E_q[ \cdot]$ are computed with respect to the marginal distribution $q(x(t))$ of the approximating process.


The KL-divergence is minimised using the mean and covariance differential equations (\ref{eq:vb_mean}) and (\ref{eq:vb_cov}) as constraints. 
Introducing Lagrange multiplier functions $\lambda(t)$ and $\Psi(t)$ for the constraints, the objective function is given by 
\begin{align}
\mathcal{F}(A,b,m,P) & = \int_{t_0}^{t_K} \left\{\E_q\left[e(x,t) +\sum_{k=1}^K u_k(x) \delta(t-t_k)\right] \right. \nonumber \\
 & \left. - \lambda(t)^T\left(\frac{dm}{dt}+A(t)m(t) - b(t)\right) \right.\nonumber \\
 & \left. - \trace\left[ \Psi(t) \left(\frac{dP}{dt}+A(t)P(t)+P(t)A(t)^T-\Sigma(t)\right) \right] \right\}\, dt.
\label{eq:KL_div3}
\end{align}
For notational convenience, we use $x = x(t)$. 
The Euler--Lagrange equations for the optimal $A(t)$ and $b(t)$ can then be written as [7,8]
\begin{align}
\frac{d}{dt}\lambda(t) & = A^T(t)\lambda(t) -\nabla_m \E_q[e(x,t)]] \label{eq:dlambda}\\
\frac{d}{dt} \Psi(t)& =  \Psi(t) A(t) + A^T(t) \Psi(t)- \nabla_P \E_q[e(x,t)] \label{eq:dPsi} \\
A(t) & = - \E_q[F_x(x,t)] + 2 \Sigma(t) \Psi(t) \label{eq:dA}\\
b(t) & = \E_q[f(x,t)] + A(t)m(t) - \Sigma(t) \lambda(t) \label{eq:db}.
\end{align}

At observation times, the Lagrange multiplier functions satisfy 
\begin{align}
\lambda(t_k^+) & = \lambda(t_k^-) + \nabla_m\E_q[u_k(x)], \label{eq:bc_lambda}\\
\Psi(t_k^+) & = \Psi(t_k^-) + \nabla_P \E_q[u_k(x)]. \label{eq:bc_psi}
\end{align}

To find a solution satisfying the Euler--Lagrange equations, we propose here a slight modification of the iterative algorithm given by Archambeau \emph{et al.} \cite{Archambeau2007, Archambeau2008}. 
Given the previous estimates $A^{(k)}(t)$ and $b^{(k)}(t)$, the mean and covariance differential equations (\ref{eq:vb_mean}) and (\ref{eq:vb_cov}) are solved forward in time from $t_0$ which gives $m^{(k+1)}(t)$ and $P^{(k+1)}(t)$. 
Using these, the Lagrange differential equations (\ref{eq:dlambda}) and (\ref{eq:dPsi}) are then solved backward in time from $t_K$ which gives $\lambda^{(k+1)}(t)$ and $\Psi^{(k+1)}(t)$. 
New estimates $A^{(k+1)}(t)$ and $b^{(k+1)}(t)$ are then computed by using a damped fixed-point update 
\begin{align}
A^{(k+1)}(t) & = A^{(k)}(t) + \gamma_k \left(A(t) - A^{(k)}(t) \right)\\
b^{(k+1)}(t) & = b^{(k)}(t) + \gamma_k \left(b(t) - b^{(k)}(t)\right),
\end{align}
where $A(t)$ and $b(t)$ are computed using Equations (\ref{eq:dA}) and (\ref{eq:db}) respectively, and $\gamma_k \in (0, 1)$ is a damping parameter that is used to prevent numerical instabilities caused by too large updates. 
Instead of using constant damping parameter as in \cite{Archambeau2007}, we propose to select the parameter $\gamma_k$ using backtracking line search so that the objective function and therefore also the KL-divergence is reduced at each time step. 
In a computer implementation of the variational Gaussian smoother, the values of the functions $A(t)$ and $b(t)$ are computed and stored at discrete time points. 
The differential equations can be solved using any standard numerical solver such as Euler or Runge--Kutta methods. 

\subsection{Variational Gaussian smoothing for singular models}

Note that to compute the term $e(x,t)$, a nonsingular effective diffusion matrix $\Sigma(t)$ is required. 
We extend here the variational Gaussian smoother to singular models of the form
\begin{align}
\frac{d x_1}{dt} & = F_1(t)x \label{eq:singular_f1} \\
d x_2 & = f_2(x, t) \, dt + \sqrt{\Sigma_2(t)} \, d \beta, \label{eq:singular_f2}
\end{align}
where $\Sigma_2(t)$ is a nonsingular diffusion matrix.

Denote by $n_1$ and $n_2$ the dimensions of the state vectors $x_1$ and $x_2$ respectively. 
We now seek an approximate smoothing process of the form
\begin{align}
\frac{d x_1}{dt} & = F_1(t)x \\
d x_2 & = (-A(t)x  + b(t))\, dt + \sqrt{\Sigma_2(t)} \, d \beta,
\end{align}
where $A(t)$ and $b(t)$ are the variational parameters.
This gives a Gaussian process, where mean $m(t)$ and covariance $P(t)$ follow the differential equations
\begin{align}
\frac{d}{dt}m(t) & = -\tilde{A}(t)m(t)+\tilde{b}(t) \label{eq:vb_mean2}\\
\frac{d}{dt}P(t) & = -\tilde{A}(t)P(t) - P(t)\tilde{A}^T(t) + \Sigma(t), \label{eq:vb_cov2}
\end{align}
with
\begin{equation}
\tilde{A}(t) = \left[ \begin{array}{c}  -F_1(t)\\ A(t) \end{array} \right], \quad \tilde{b}(t) = \left[ \begin{array}{c} 0_{n_1 \times 1} \\ b(t) \end{array} \right], \quad \Sigma(t) = \left[ \begin{array}{cc} 0_{n_1 \times n_1} & 0_{n_1 \times n_2} \\ 0_{n_2 \times n_1} & \Sigma_2(t) \end{array} \right].
\end{equation}
The KL-divergence term for the singular system can be computed using Girsanov's theorem and is given by (see \cite{SarkkaSottinen2007} and Appendix A for details) 
\begin{equation}
\KL(\mathbb{Q}_X \, || \,  \mathbb{P}_{X\cond Y}) = \int_{t_0}^{t_K}  \E_q\left[e(x, t) + \sum_{k=1}^K u_k(x)\delta(t-t_k) \right] \, dt,
\label{eq:KL_div_singular}
\end{equation}
where $u_k(x)$ is given by eq. (\ref{eq:uk}) and
\begin{align}
e(x, t) & = \frac{1}{2}[f_2(x,t)+A(t)x(t) -b(t)]^T \Sigma_2^{-1}(t)[f_2(x, t)+A(t)x(t) -b(t)] .
\end{align}
The objective function is given by
\begin{align}
\mathcal{F}(A, b, m, P) & = \int_{t_0}^{t_K} \left\{\E_q\left[e(x,t) +\sum_{k=1}^K u_k(x)\delta(t-t_k)\right] \right. \nonumber \\
 & - \left. \lambda(t)^T\left(\frac{dm}{dt}+\tilde{A}(t)m(t) - \tilde{b}(t)\right) \right. \nonumber \\
 & \left. - \trace\left[ \Psi(t) \left(\frac{dP}{dt}+\tilde{A}(t)P(t)+P(t)\tilde{A}(t)^T-\Sigma(t)\right) \right] \right\}\, dt
\end{align}
and the solution satisfies
\begin{align}
\frac{d}{dt}\lambda(t) & = \tilde{A}^T(t)\lambda(t) - \nabla_m \E_q[e(x,t)]] \label{eq:dlambda2}\\
\frac{d}{dt} \Psi(t)& =  \Psi(t) \tilde{A}(t) + \tilde{A}^T(t) \Psi(t)- \nabla_P \E_q[e(x,t)] \label{eq:dPsi2} \\
A(t) & = - \E_q[F_{2,x}(x,t)] + 2 \Sigma_2(t) M\Psi(t) \\
b(t) & = \E_q[f_2(x,t)] + A(t)m(t) - \Sigma_2(t) M\lambda(t),
\end{align}
where $F_{2,x}(x,t)$ is the Jacobian of $f_2(x,t)$ and the matrix $M$, which selects the relevant part of the Lagrange parameter functions, is given by
\begin{equation}
M = \left[ \begin{array}{cc} 0_{n_1 \times n_2} & I_{n_2 \times n_2} \end{array} \right].
\end{equation}
At observation times, the Lagrange multipliers satisfy the boundary conditions given by eqs. (\ref{eq:bc_lambda})-(\ref{eq:bc_psi}).

\subsection{Gaussian filtering based Gaussian smoother}

The Gaussian filtering based Gaussian smoother uses the Gaussian approximation for the filtering distribution to form the Gaussian approximation for the smoothing distribution \cite{SarkkaSarmavuori2013}. 

The Gaussian filtering (or Gaussian assumed density filtering) approach is well known in the literature (see e.g. \cite{Kushner1967, Jazwinski1970, ItoXiong2000}) and uses the approximation 
\begin{equation}
p(x(t) \cond y_{1}, \ldots, y_{k}) \approx \text{N}(x(t) \cond m_f(t), P_f(t)), 
\label{eq:filter}
\end{equation}
where $p(x(t) \, | \, y_{1}, \ldots, y_{k})$ is the filtering distribution. 
The mean and covariance function are recursively computed using the following prediction and update steps. 
In the prediction step the mean and covariance functions are propagated from time $t_{k-1}$ to time $t_k$ using 
\begin{align}
\frac{dm_f}{dt} & = \E_f[f(x,t)] \label{eq:dmf}  \\
\frac{dP_f}{dt} & = \E_f[(x-m_f)f(x,t)^T] + \E_f[f(x,t)(x-m_f)^T] + \Sigma(t). \label{eq:dPf}
\end{align}
In the update step, the information from the latest measurement $y_k$ is used to update the predicted estimates $m(t_k^-)$ and $P(t_k^-)$ using equations
\begin{align}
S_k  &= \E_f\left[(h_k(x) - \E_f[h_k(x)] )(h_k(x)-\E_f[h_k(x)] )^T\right] + R_k  \\
K_k & = \E_f\left[(x-m_f(t_k^-))(h_k(x)-\E_f[h_k(x)] )^T\right] S^{-1}_k \\
m_f(t_k) & = m_f(t_k^-) + K_k\left(y_k - \E_f[h_k(x)]\right) \\
P_f(t_k) & = P_f(t_k^-) - K_kS_k K_k^T.
\end{align}

S\"arkk\"a and Sarmavuori \cite{SarkkaSarmavuori2013} extend the general Gaussian filtering ideas also to smoothing problems and derive three types of Gaussian smoothers labeled as type I, type II and type III. 
The type I smoother is derived by using the Gaussian approximation for the filtering distribution to approximate the exact partial differential equations for the smoothed mean and covariance. 
The type II smoother is derived by discretizing the dynamic model, applying the discrete-time smoothing equations and then taking the limit as the discretization time approaches zero. 
Formulating the type II smoothing equations into a computationally efficient form gives the type III smoother. 

The type I smoothing equations are numerically quite sensitive, which causes the type I smoother to diverge quite often compared to the type II and type III smoothers.
Also, the type I smoother is computationally more demanding and was not found to be clearly better than the type II and type III smoothers in the synthetic-data example considered in \cite{SarkkaSarmavuori2013}. 
For this reason we decided to concentrate on the type II and type III smoothers in this paper. 

The type II smoothing equations are given by
\begin{align}
\frac{dm_s}{dt} & = \E_f[f(x,t)] + \left[\E_f[f(x,t)(x-m_f)^T] + \Sigma(t) \right]P_f^{-1}(m_s-m_f) \label{eq:dms} \\
\frac{dP_s}{dt} & = (\E_f[f(x,t)(x-m_f)^T] + \Sigma)P_f^{-1}P_s \nonumber \\
& + P_sP_f^{-1}(\E_f[f(x,t)(x-m)^T]^T+ \Sigma(t)) - \Sigma(t) \label{eq:dPs}.
\end{align} 
If the Jacobian $F_x(x,t)$ of $f(x,t)$ is available, we can alternatively use
 \begin{equation}
\E_f[f(x,t)(x-m_f)^T] =  \E_f[F_x(x,t)]P_f  \label{eq:type2_jacobian}
\end{equation}
in the filtering and smoothing equations. 
The expectations in the smoothing equations are with respect to the filtering density, which means that they can be computed already during the filtering stage.
This is used in the type III smoothing equations, which reformulate the filtering and type II smoothing equations so that no expectations need to be computed during the smoothing stage (see \cite{SarkkaSarmavuori2013} for details).

\subsection{Differences between the Gaussian smoothers}

In this section, we study the differences between the variational and Gaussian filtering based Gaussian smoothers.
First we introduce a change of variables that converts the type II Gaussian smoothing equations to a form similar to the variational Gaussian smoothing equations. 
This conversion is similar to the results from Rauch, Tung and Striebel \cite{RauchTungStriebel1965}, where it was shown that the Rauch--Tung--Striebel smoothing equations are formally equivalent to the smoothing equations presented by Bryson and Frazier \cite{BrysonFrazier1963}.

The conversion is achieved by using the change of variables
\begin{equation}
\lambda(t) = -P_f^{-1}(t)(m_s(t)-m_f(t)), \, \Psi(t) = -\frac{1}{2}(P_f^{-1}(t) - P_s^{-1}(t)) \label{eq:type2_lagrange}. 
\end{equation}
Inserting $\lambda(t)$ and $\Psi(t)$ to the Equations (\ref{eq:dms})-(\ref{eq:dPs}) and computing the time derivatives of $\lambda(t)$ and $\Psi(t)$ gives (see Appendix B)
\begin{align}
\frac{d}{dt}m_s(t) & = -A(t)m_s(t)+b(t) \label{eq:type2_mean2} \\
\frac{d}{dt}P_s(t) & = -A(t)P(t) - P(t)A^T(t) + \Sigma(t) \label{eq:type_cov2} \\
\frac{d}{dt}\lambda(t) &= A^T(t) \lambda(t) - 2 \Psi(t) \Sigma(t) \lambda(t) \label{eq:type2_lambda} \\
\frac{d}{dt}\Psi(t) &  = \Psi(t)A(t) + A^T(t)\Psi(t) - 2\Psi(t) \Sigma(t) \Psi(t), \label{eq:type2_psi}
\end{align}
where 
\begin{align}
A(t) & = -\E_f[F_x(x,t)] + 2 \Sigma(t) \Psi(t) \label{eq:type2_A} \\
b(t) & = \E_f[f(x,t)] + \E_f[F_x(x,t)](m_s(t)-m_f(t)) + A(t) m_s(t) - \Sigma(t) \Psi(t) \label{eq:type2_b} .
\end{align}
For linear measurement function $h_k(x) = H_kx$, the measurement update for $\lambda(t)$ and $\Psi(t)$ in the variational form of the type II smoother is the same as in the variational Gaussian smoother and is given by
\begin{align}
\lambda(t_k^+) & = \lambda(t_k^-) + H_k^T R_k^{-1}(m_s(t_k)-y_k) \\ 
\Psi(t_k^+) & = \Psi(t_k^-) + \frac{1}{2} H_k^T R_k^{-1} H_k. 
\end{align}
For nonlinear measurement function, the measurement updates for $\lambda(t)$ and $\Psi(t)$ are in general not equal to the variational Gaussian smoother update Equations (\ref{eq:bc_lambda}) and (\ref{eq:bc_psi}). 

Similarities to the variational Gaussian smoothing equations are evident from Equations (\ref{eq:type2_mean2})-(\ref{eq:type2_b}). 
Note that Equations (\ref{eq:type2_A}) and (\ref{eq:type2_b}) for the parameters $A(t)$ and $b(t)$ are otherwise similar to the variational Gaussian smoothing Equations (\ref{eq:dA}) and (\ref{eq:db}), but the function $f(x,t)$ is replaced with statistical linearization with respect to the filtering distribution:
\begin{equation}
f(x,t) \approx \E_f[f(x,t)] + \E_f[F_x(x,t)](x-m_f(t)) \label{eq:stat_lin}. 
\end{equation}
Furthermore, using the statistical linearization (\ref{eq:stat_lin}) to approximate the gradients in Equations (\ref{eq:dlambda}) and (\ref{eq:dPsi}) gives
\begin{align}
\nabla_m \E[e(x,t)] & \approx 2 \Psi(t) \Sigma(t) \lambda(t) \\
\nabla_P \E[e(x,t)] & \approx 2 \Psi(t) \Sigma(t) \Psi(t). 
\end{align}
That is, the differential equations (\ref{eq:type2_lambda}) and (\ref{eq:type2_psi}) can be seen as an approximation to the exact differential equations (\ref{eq:dlambda}) and (\ref{eq:dPsi}) in the variational Gaussian smoother. 
This suggests that for linear measurements, the type II Gaussian smoother can be seen to approximate the variational Gaussian smoother by using statistical linearization with respect to the filtering distribution. 

We can use the type II Gaussian smoothing solution as an initial iterand for the variational parameters by computing $A^{(0)}(t)$ and $b^{(0)}(t)$ using equations (\ref{eq:type2_A}) and (\ref{eq:type2_b}), where $\lambda(t)$ and $\Psi(t)$ are given by (\ref{eq:type2_lagrange}). 
In \cite{Vrettas2010} it was noted that the iterative solution of the variational Gaussian smoothing equations is sensitive to the initial values of the variational parameters $A(t)$ and $b(t)$. 
This way, the variational Gaussian smoother can be seen as an iterative way to improve the type II Gaussian smoothing solution.
The benefit of using the variational Gaussian smoother to improve the type II smoother results is studied further in the synthetic-data examples.

\section{Computation of Gaussian expectations in the variational Gaussian smoother}

The Gaussian smoothers considered in this paper require computations of Gaussian expectations over arbitrary nonlinear functions.
For some simple models these can be computed analytically, but for many models no analytical expression exists.
The computation of Gaussian expectations for the Gaussian filtering based smoothers is presented in \cite{SarkkaSarmavuori2013}. 
In this work we concentrate on computing the expectations in the variational Gaussian smoothing equations and present the extended, cubature, unscented and Gauss--Hermite forms of the variational Gaussian smoother. 

Gaussian expectations need to be computed in the differential equations (\ref{eq:dlambda}) and (\ref{eq:dPsi}) for the functions $\lambda(t)$ and $\Psi(t)$, in Equations (\ref{eq:dA}) and (\ref{eq:db}) for the variational parameter functions $A(t)$ and $b(t)$ and in the measurement update equations (\ref{eq:bc_lambda}) and (\ref{eq:bc_psi}). 
In order to avoid computing derivatives of the drift function $f(x,t)$ and of the measurement function $h_k(x)$, the gradients with respect to $m(t)$ and $P(t)$ can be written in the form (see Appendix C)
\begin{align}
\nabla_m \E[e(x,t)] & =  P^{-1}\E[e(x,t)(x-m)]  \\
\nabla_P \E[e(x,t)] & = \frac{1}{2}P^{-1}\E[e(x,t)(x-m)(x-m)^T]P^{-1} - \frac{1}{2}\E[e(x,t)]P^{-1}.
\end{align}
For the measurement updates, the gradients are computed similarly with $e(x,t)$ replaced by $u_k(x)$. 

If the Jacobians $F_x(x,t)$ and $H_{k,x}(x)$ of $f(x,t)$ and $h_k(x)$ are available, the gradients can be alternatively written in the form
\begin{align}
\nabla_m \E[e(x,t)] & =  \E[e_x(x,t)] \label{eq:dm_alt} \\
\nabla_P \E[e(x,t)] & = \frac{1}{2}P^{-1}\E[e_x(x,t)(x-m)^T]^T. \label{eq:dP_alt} 
\end{align}
where 
\begin{equation}
e_x(x,t) = \nabla_x e(x,t) = [F_x(x,t) + A(t)]^T\Sigma^{-1}(t)[f(x,t) + A(t)x(t) -b(t)].
\end{equation}
To compute the gradients in the measurement update, $e_x(x,t)$ is replaced with
\begin{equation}
u_{k,x}(x) = \nabla_x u_k(x) =H_{k,x}^T(x) R_k^{-1}[h_k(x) - y_k]. 
\end{equation}

\subsection{Taylor series based linearization}
In the extended Kalman filter and smoother, the Gaussian expectations are computed by using a first order Taylor series linearization. 
Proceeding similarly, we use the following approximations for the extended variational Gaussian smoother: 
\begin{align}
f(x,t) & \approx f(m, t) + F_x(m,t)(x(t)-m(t)) \\
h_k(x) & \approx h_k(m) + H_{k,x}(m)(x(t))-m(t)),
\end{align}
where $F_x(m,t)$ and $H_{k,x}(m)$ are the Jacobians of $f(x,t)$ and $h_k(x)$ respectively evaluated at the current mean estimate $m(t)$.  
Using this approximation, the expectations in Equations (\ref{eq:dA}) and (\ref{eq:db}) are given by
\begin{align}
\E[f(x,t)] & \approx f(m,t) , \\
\E[F_x(x,t)] & \approx F_x(m,t). 
\end{align}
The expectations needed in the gradients are given by
\begin{align}
\nabla_m \E[e(x,t)] & \approx (F_x(m,t) + A(t))^T\Sigma^{-1}(t)(f(m,t) + A(t)m(t) - b(t))  \\
\nabla_P \E[e(x,t)] & \approx \frac{1}{2}(F_x(m,t) + A(t))^T\Sigma^{-1}(t)(F_x(m,t) + A(t))
\end{align}
and for the measurement update 
\begin{align}
\nabla_m \E[u_k(x)] & \approx H^T_{k,x}(m)R_k^{-1}(t)(h_k(m)  - y_k)  \\
\nabla_P \E[u_k(x)] & \approx \frac{1}{2}H^T_{k,x}(m)^TR_k^{-1}(t)H_{k,x}(m).
\end{align}

\subsection{Sigma-point methods}

The general sigma-point rule computes the Gaussian expectations using the approximation
\begin{equation}
\E[g(x,t)] \approx \sum_i W^{(i)} g(m + \sqrt{P} \xi_i,t),
\end{equation}
where the weights $W^{(i)}$ and vectors $\xi_i$ are chosen depending on the used sigma-point method. 
In this paper we consider the cubature, unscented and Gauss--Hermite sigma-point methods. 
The cubature method uses $2n$ sigma-points with vectors
\begin{equation}
\xi_i = \left\{ \begin{array}{ll} \sqrt{n} e_i, & i = 1,\ldots,n \\ -\sqrt{n} e_{i-n}, & i=n+1,\ldots,2n \end{array} \right.
\label{eq:ckf_xi}
\end{equation}
and weights $W^{(i)} = 1/(2n)$ for all $i=1,\ldots,2n$. 
The unscented transform uses $2n+1$ sigma-points with vectors
\begin{align}
\xi_0 & = 0, \quad \xi_i  = \left\{ \begin{array}{ll} \sqrt{\lambda + n} e_i, & i = 1,\ldots,n \\ -\sqrt{\lambda +n} e_{i-n}, & i=n+1,\ldots,2n.\end{array} \right. ,
\label{eq:ukf_xi}
\end{align}
where $\lambda = \alpha^2(n+\kappa) - n$ and $\alpha$, $\beta$ and $\kappa$ are parameters of the method. 
The weights are defined to be
\begin{align}
W_m^{(0)} & = \frac{\lambda}{n + \lambda}, \quad W_m^{(i)} = \frac{1}{2(n + \lambda)}, \,\, i = 1,\ldots,2n \\
W_c^{(0)} & = W_m^{(0)} + 1-\alpha^2 + \beta, \quad W_c^{(i)} = W_m^{(i)}, \,\, i=1,\ldots,2n.
\end{align}
The weights $W_m$ are used in approximating the transformed mean and $W_c$ in the covariance approximation.
The cubature method is a special case of the unscented transform with parameters $\alpha=1$, $\beta = \kappa =0$. 
The Gauss--Hermite integration uses $s^n$ sigma-points, where $s$ is a parameter that gives the order of the used Hermite polynomial. 
Details for computing the vectors $\xi_i$ and weights $W^{(i)}$ are given in \cite{ItoXiong2000, WuEtAl2006}. 

The sigma-point approximation for expectations in Equations (\ref{eq:dA}) and (\ref{eq:db}) is given by
\begin{align}
\E[f(x,t)] & \approx \sum_i W^{(i)} f(m + \sqrt{P}\xi_i, t) \label{eq:Ef_approx} \\
\E[F_x(x,t)] &= \E[f(x,t)(x-m)^T]P^{-1}  \approx \sum_i W^{(i)} f(m + \sqrt{P}\xi_i,t)\xi_i^T \sqrt{P}^{-1} \label{eq:EF_approx}.
\end{align}
The general sigma-point approximations for the gradients are given by
\begin{align}
\nabla_m \E[e(x,t)] & \approx \sum W^{(i)} e(m+\sqrt{P}\xi_i, t)\sqrt{P}^{-1}\xi_i   \\
\nabla_P \E[e(x,t)] & \approx \frac{1}{2}\sum W^{(i)} e(m + \sqrt{P} \xi_i,t)\sqrt{P}^{-T}\left( \xi_i \xi_i^T -I \right) \sqrt{P}^{-1}.
\end{align}
The expectations needed in the observation updates (\ref{eq:bc_lambda}) and (\ref{eq:bc_psi}) are computed with $e(x,t)$ replaced by $u_k(x)$. 

The sigma-point approximation for the alternative forms (\ref{eq:dm_alt}) and (\ref{eq:dP_alt}) of the gradients are given by
\begin{align}
\nabla_m \E[e(x,t)] & \approx \sum W^{(i)} e_x(m+\sqrt{P}\xi_i), \label{eq:dm_alt2} \\
\nabla_P \E[e(x,t)] & \approx \frac{1}{2}\sum W^{(i)} \sqrt{P}^{-T}\xi_i e_x^T(m + \sqrt{P} \xi_i,t). \label{eq:dP_alt2}
\end{align}
The measurement updates are computed similarly, with $e_x(x,t)$ replaced by $u_{k,x}(x)$. 
The cubature, unscented and Gauss--Hermite forms of the variational Gaussian smoother are then obtained by using the corresponding choice for the weights $W^{(i)}$ and vectors $\xi_i$. 

For a linear drift function $f(x,t)$, the term $\E[e(x,t)(x-m)(x-m)^T]$ is a fourth order polynomial. 
For this reason, the sigma-point rules that are only accurate up to a third order monomial (cubature and unscented rule) give generally a poor approximation of this expectation.
Therefore, for cubature and unscented sigma-point methods, the use of the alternative form given by Equations (\ref{eq:dm_alt2}) and (\ref{eq:dP_alt2}) is recommended. 


\section{Numerical experiments}

The Gaussian smoothers are compared using two different synthetic-data experiments. 
The tests are done by running first the Gaussian filtering based Gaussian smoother (GFGS) and then using the result as initial conditions for the variational Gaussian smoother (VGS). 
The VGS iteration is terminated when the absolute change in the KL-divergence between successive iterations is less than $10^{-3}$. 
 
In the first experiment, a one dimensional double well system is used. 
The same system was also used to demonstrate the VGS in \cite{Archambeau2007, Shen2010, Vrettas2010}.
A 5-dimensional reentry problem is used for the second experiment. 
This system was used to test the continuous-discrete unscented Kalman filter in \cite{Julier2004} and demonstrates the use of VGS for singular systems. 
Also, for this system the computation of the needed Gaussian expectations is not possible analytically.

The metrics used to compare the estimates given by the Gaussian smoothers are the root mean square error (RMSE), negative log-likelihood (NLL) and 95\%-consistency.
The RMSE and NLL are given by equations
\begin{align}
\text{RMSE} &  = \sqrt{\frac{1}{t_K-t_0}\int_{t_0}^{t_K} \| x(t)-m(t) \|^2 \, dt} \\
\text{NLL} & = \frac{1}{t_K-t_0}\int_{t_0}^{t_K} \ln \text{N}(x(t) \, | \, m(t), P(t)) \, dt,
\end{align}
where $x(t)$ is the true state, $m(t)$ is the estimated mean and $P(t)$ is the estimated covariance. 
The 95\%-consistency is defined as the fraction of times the true state is inside the 95\% ellipsoid of $\text{N}(m(t), P(t))$. 
The values of the continuous time metrics are computed using the values of $m(t)$ and $P(t)$ computed at discrete time points. 

The different approximation methods for computing the Gaussian expectations in the smoothing equations are also compared.
The tested methods are labeled as
\begin{itemize}
\item EXT: The method using the Taylor series based linearization.
\item CT: Cubature rule based sigma-point method. 
\item UT: Unscented rule based sigma-point method with paremeter values $\alpha = 1$, $\beta = 2$ and $\kappa = 0$.
\item G-H: Gauss--Hermite series based sigma-point method with order 3. 
\end{itemize}
The methods labeled CT2, UT2 and G-H2 use the respective sigma-point method with the alternative formulation given in Equations (\ref{eq:dm_alt2}) and (\ref{eq:dP_alt2}). 


\subsection{Double well}
The double well system is given by 
\begin{align}
dx & = 4x(1-x^2) dt + \sqrt{\sigma} d\beta, \quad y_k = x(t_k) + v_k,
\label{eq:dw}
\end{align}
where the measurement noise $v_k$ is zero-mean Gaussian with variance $R$. 
The prior distribution for the double well system is non-Gaussian and multimodal, but the smoothing distribution can be reasonably well approximated with a Gaussian provided that the measurement variance is not too large. 
The modes are located at $x=1$ and $x=-1$ and for sufficiently large value of the process noise parameter $\sigma$, there is frequent transition from one mode to the other. 

The Gaussian smoothers are compared for 4 different values for the measurement variance $R$. 
For each value of the measurement variance $R$, a data set of 100 Monte Carlo simulations is generated using Euler--Maruyama discretization with time step $\Delta t = 0.01$ from $t_0 = 0$ to $t = 10$. 
The process noise parameter is chosen to be $\sigma = 1$, which is sufficiently large to cause frequent transition between the modes. 
The initial state is chosen as $x(0) \sim \text{N}(0, 1)$. 

For this system, the Gaussian expectations needed in the Gaussian smoothers can be computed analytically (see e.g. \cite{Archambeau2007}).
For both Gaussian smoothers, the differential equations are solved using 4th order Runge--Kutta method with time step $0.01$. 
Average number of 35 iterations was observed for the VGS when using the GFGS results as initial conditions. 
Also, we observed a much faster convergence using this initialization than with the naive initialization using just the initial conditions. 

\begin{figure}[h]
\centering
\includegraphics[scale=0.25]{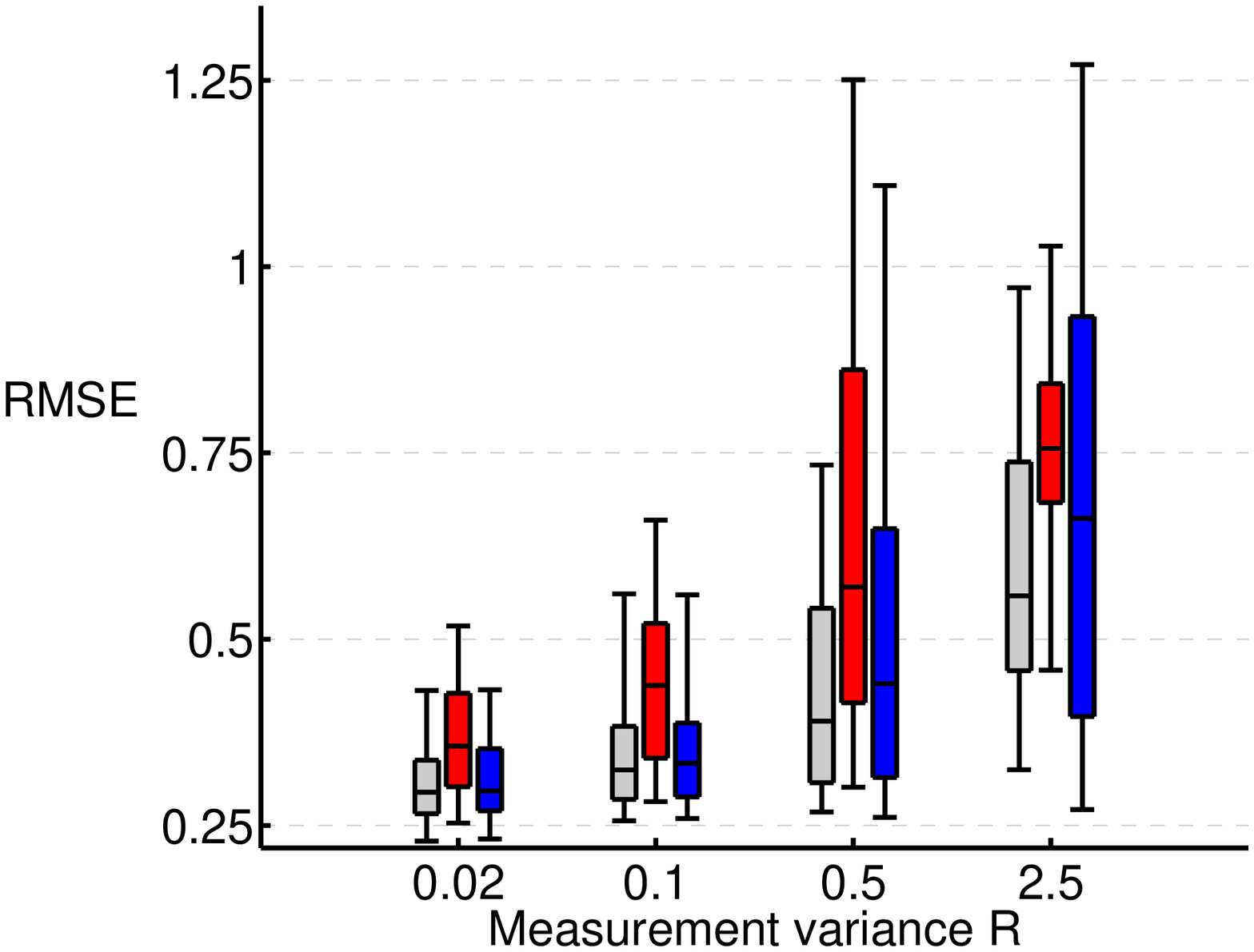} \, 
\includegraphics[scale=0.25]{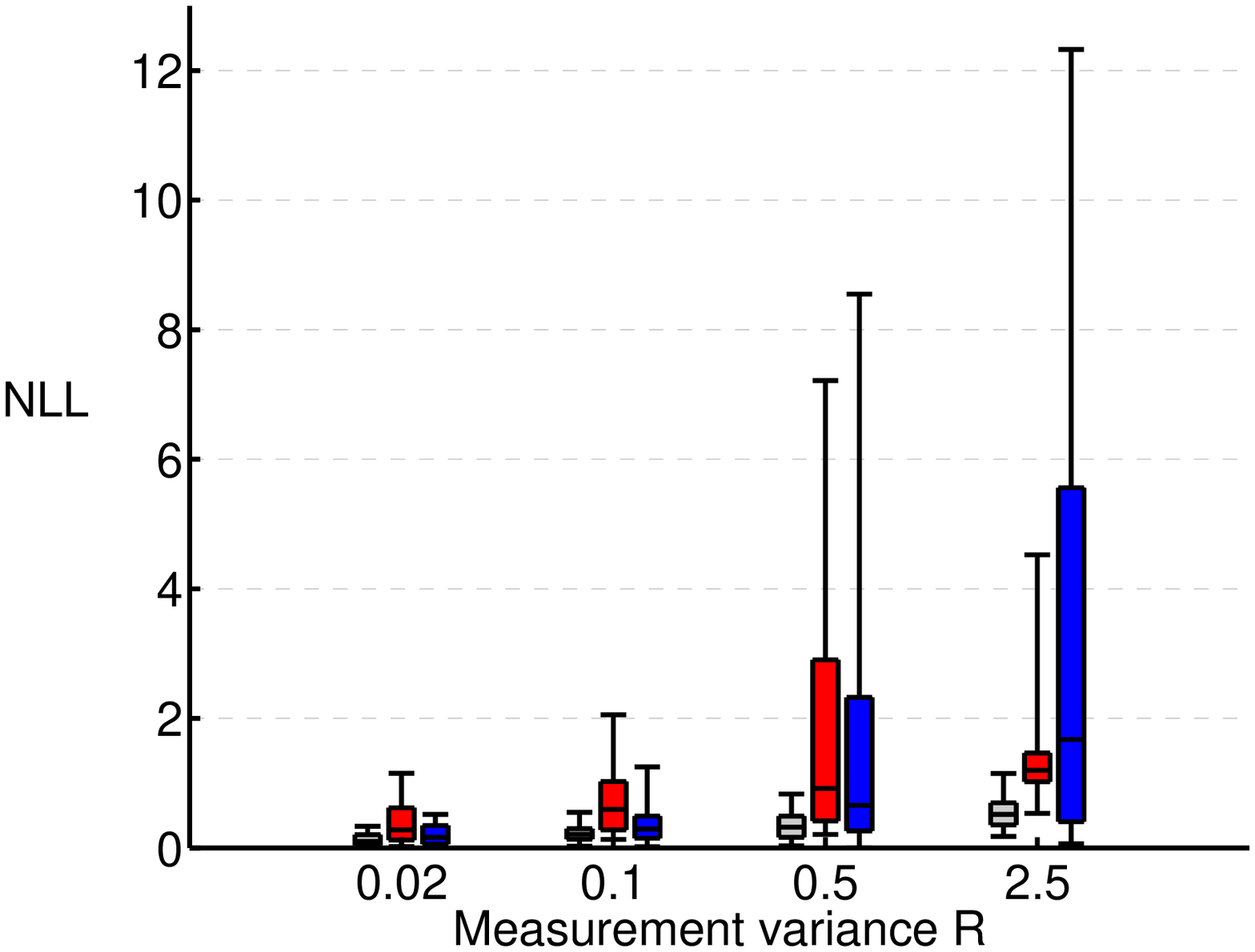}
\caption{The RMSE (left) and NLL (right) for the reference (grey), GFGS (red) and VGS (blue) for different values of the measurement variance $R$. The Gaussian expectations in the smoothing equations are computed analytically. The boxplots show the 5\%, 25\%, 50\%, 75\% and 95\% quantiles. }
\label{fig:RMSE_NLL1} 
\end{figure}

The boxplots of the RMSE and NLL results for different values of the measurement variance $R$ are shown in Figure \ref{fig:RMSE_NLL1}. 
For comparison, a reference smoothing solution is also computed using a finite difference approximation of the exact Bayesian smoothing equations \cite{Sarkka2006}. 
For the reference solution, the NLL is computed using a finite difference approximation for the smoothing density. 

\begin{figure}[h]
\centering
\includegraphics[scale=0.25]{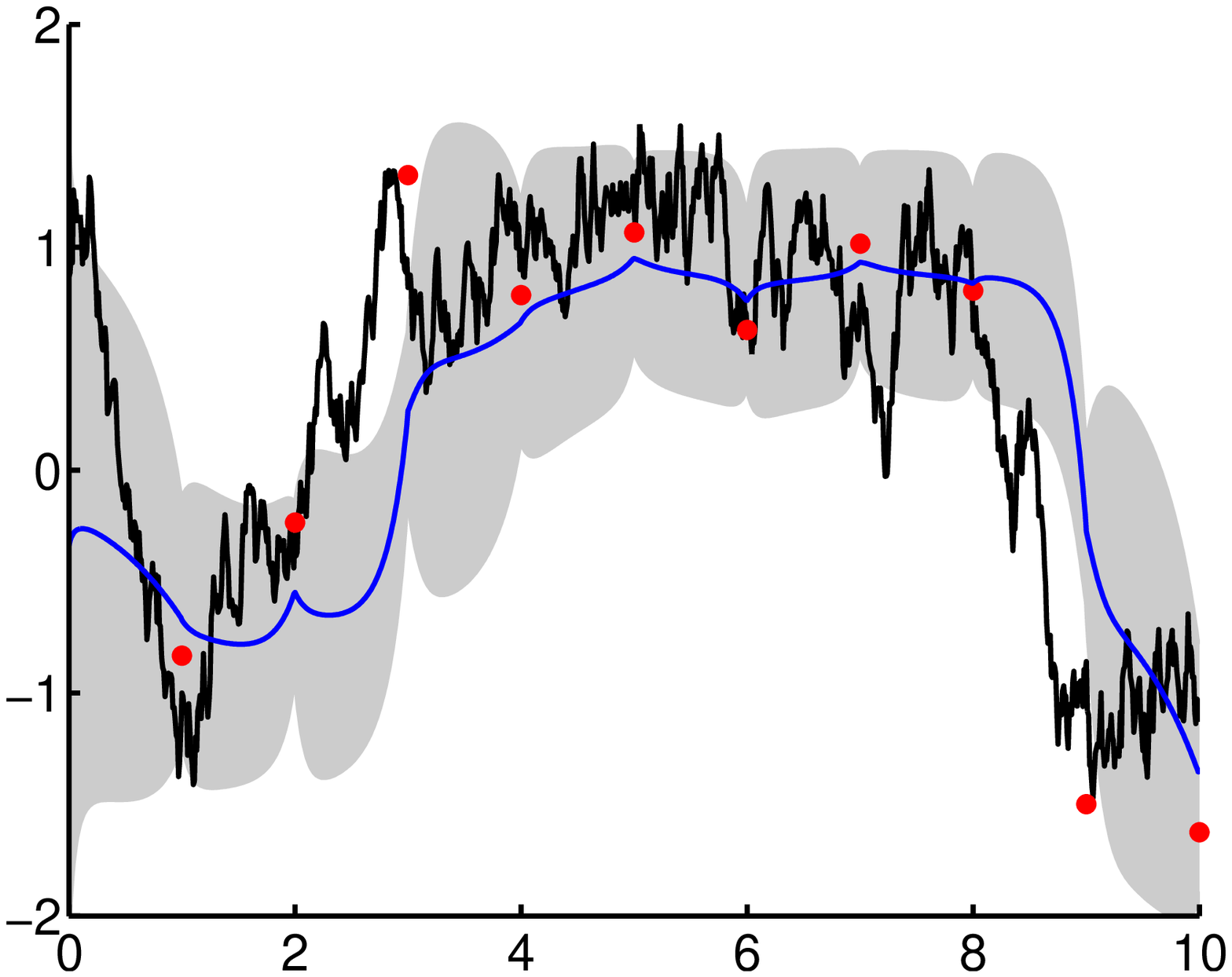} \, 
\includegraphics[scale=0.25]{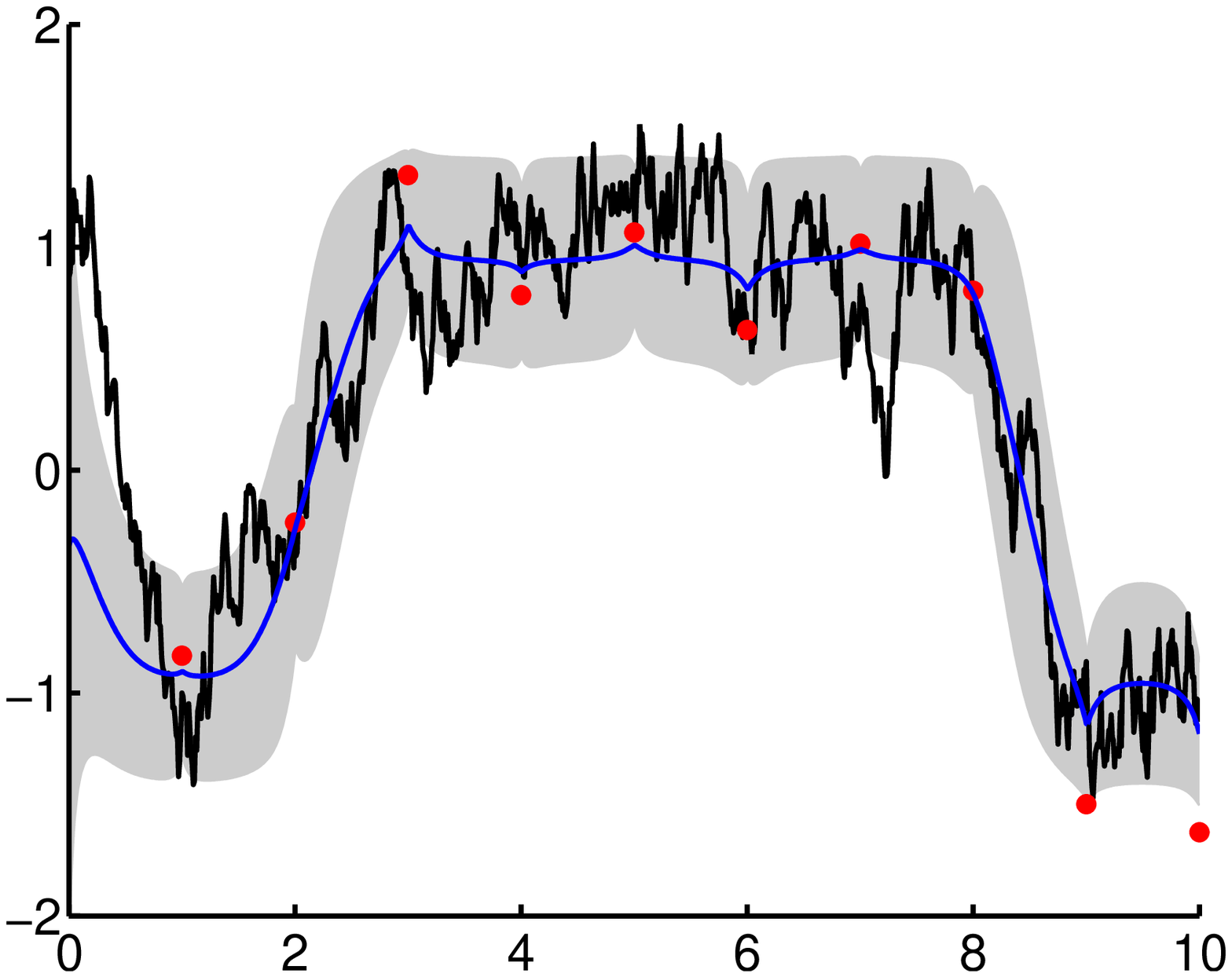}
\caption{The estimated mean (blue) and 95\%-confidence region (shaded) for GFGS (left) and VGS (right). True state (black) and measurements (red dot) are also shown. The Gaussian expectations in the smoothing equations are computed analytically. }
\label{fig:ex1} 
\end{figure}

For small values of the measurement variance, the VGS results are very close to the reference solution and clearly outperform the GFGS approach in terms of RMSE and NLL.
The relatively poor RMSE values for the GFGS are due to the poor estimation of the transitions between the two modes. 
A typical time series for measurement variance $R=0.1$ is shown in Figure \ref{fig:ex1}. 

For measurement variance $R=2.5$, the true posterior is bimodal and the VGS tends to have the estimated mean close to one of the modes with relatively small variance.
This results in very large NLL values for the VGS, when the true path is not close to the mode. 
In comparison, the GFGS tends to have the mean close to zero with the 95\%-confidence region covering both modes. 
This shows as very good NLL values and smaller spread of the RMSE values compared to the VGS. 

The mean 95\%-consistency results over the 100 Monte Carlo simulations are shown in Table \ref{table:95cons1}. 
From the consistency results, we see that in general the VGS tends to underestimate the variance compared to the GFGS. 
This is especially clear in the $R=2.5$ case.
The underestimation of the variance is a general property of the variational type of approximations \cite[p. 431]{MacKay2003}. 

\begin{table}
\caption{The mean 95\%-consistencies for the GFGS and VGS for different values of measurement variance. The Gaussian expectations in the smoothing equations are computed analytically.}
\label{table:95cons1}
\center
\begin{tabular}{r | c c c c}
& \multicolumn{4}{c}{Measurement variance $R$} \\ 
 & $0.02$ & $0.1$ &  $0.5$ &  $2.5$  \\
 \hline
GFGS & 0.93 & 0.89 & 0.81 & 0.91 \\
VGS & 0.91 & 0.88 & 0.80 & 0.74
\end{tabular}
\end{table}

The different approximation methods for the Gaussian expectations were compared for measurement variance $R=0.02$ and  using the same data set of 100 Monte Carlo simulations as in the first experiment. 
The VGS failed to converge to a solution on 5 cases using the EXT method, but no failures were observed using the other methods. 

\begin{figure}[h]
\centering
\includegraphics[scale =0.25]{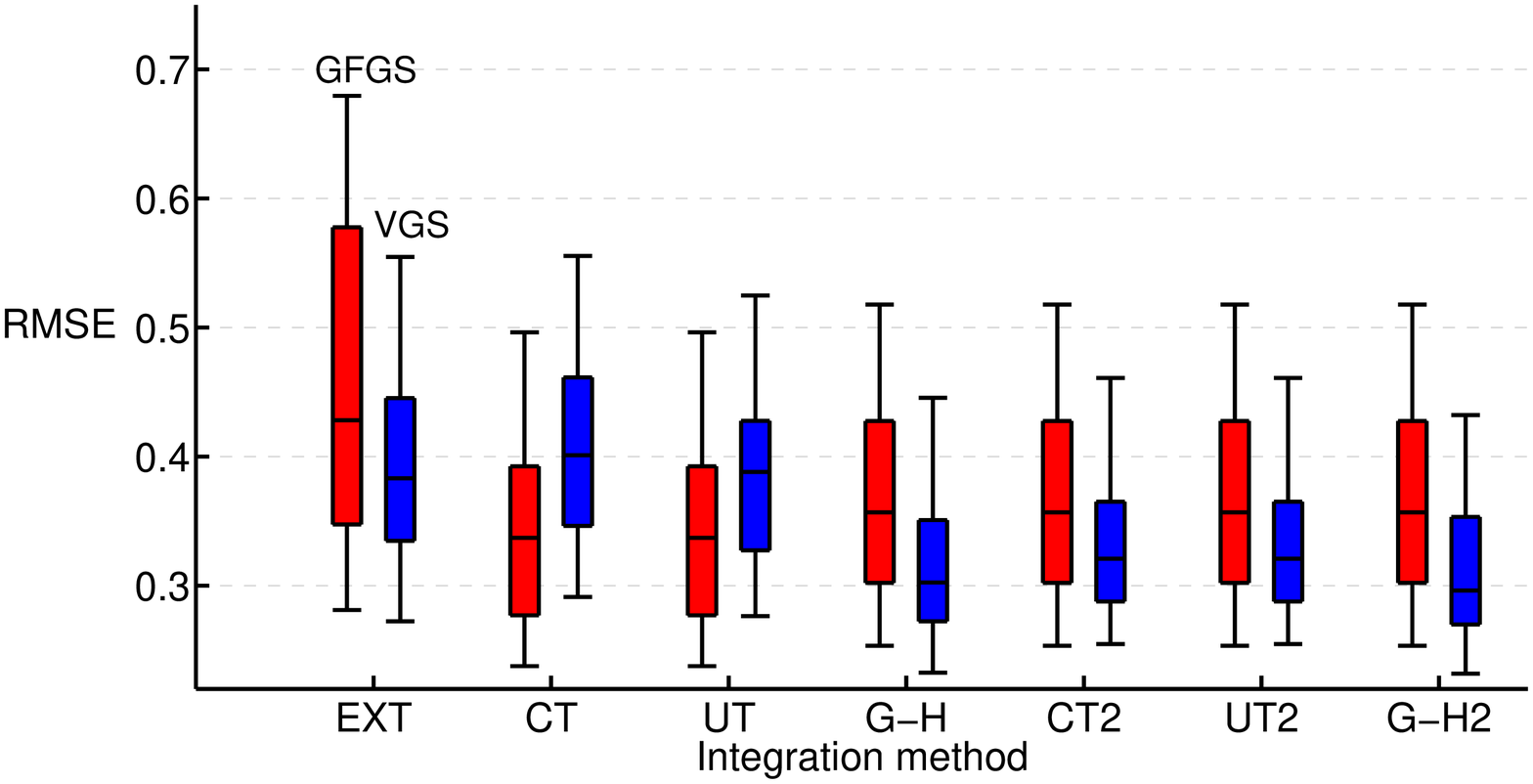}
\includegraphics[scale=0.25]{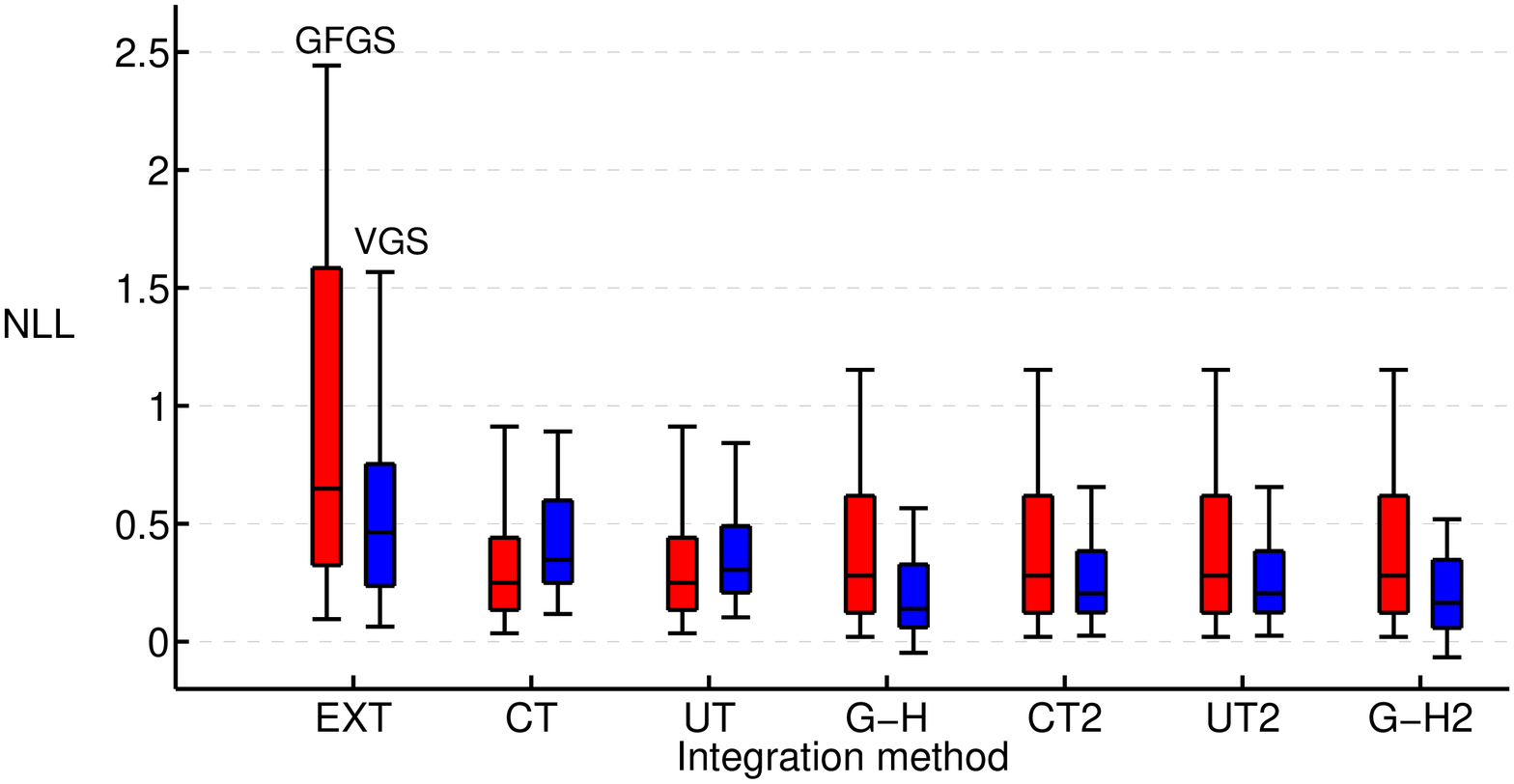} 
\caption{The RMSE (top) and NLL (bottom) for the GFGS (red) and VGS (blue) using different methods to compute the Gaussian expectations. The boxplots show the 5\%, 25\%, 50\%, 75\% and 95\% quantiles. }
\label{fig:RMSE_NLL3} 
\end{figure}

The boxplots of the RMSE and NLL values for the GFGS and VGS when using the different methods to approximate the Gaussian expectations are shown in Figure \ref{fig:RMSE_NLL3}. 
The VGS using the EXT and G-H methods clearly improve the results of the corresponding GFGS in terms of RMSE and NLL.
Numerical problems were observed for the VGS using CT and UT rules, which resulted in poor RMSE and NLL values compared to the GFGS results. 
Using the alternative formulation of CT2 and UT2 works clearly better and slightly improve the results of the GFGS using the CT and UT methods. 
The VGS using G-H2 method gives results nearly identical to the VGS using the exact Gaussian expectations. 

The mean 95\%-consistency values over the 100 Monte Carlo simulations are shown in Table \ref{table:95cons3}. 
The VGS using the EXT method shows a slight improvement in the consistency compared to the corresponding GFGS result. 
For the sigma-point methods, the differences are smaller with GFGS giving in general slightly better consistency results.  

For comparison, we also included CT2, UT2 and G-H2 versions for the GFGS that use Equation (\ref{eq:type2_jacobian}) to compute the Gaussian expectations. 
There seems to be no significant improvement in using the alternative formulation for the GFGS.
Also, this increases the computational load, since also the Jacobian needs to be evaluated for each sigma-point.



\begin{table}
\caption{The mean 95\%-consistencies for the GFGS and VGS using different methods to compute the Gaussian expectations.}
\label{table:95cons3}
\center
\begin{tabular}{r | c c c c c c c}
& \multicolumn{5}{c}{Integration method} \\ 
 &  EXT & CT  & UT & G-H &CT2 & UT2&  G-H2  \\
 \hline
GFGS & 0.83 & 0.94 & 0.94 & 0.93 & 0.93 & 0.93 & 0.93 \\ 
VGS & 0.86 & 0.93 & 0.94 & 0.93 & 0.94 & 0.94 & 0.91 
\end{tabular}
\end{table}

\subsection{Reentry}

The state $x=[r, v, \alpha]^T$ of the reentry problem consists of the vehicle's position $r$ and velocity $v$ in a 2-dimensional coordinate system and a parameter $\alpha$ of its aerodynamic properties. 
The dynamics are given by
\begin{align}
d\left[ \begin{array}{c} r(t) \\ v(t) \\ \alpha(t) \end{array} \right]  = \left[ \begin{array}{ccc} 
0_{2 \times 2} & I_{2 \times 2} & 0_{1 \times 1} \\
G(x,t) I_{2 \times 2} & D(x,t) I_{2 \times 2} & 0_{1 \times 1} \\
0_{1 \times 2} & 0_{1 \times 2} & 0_{1 \times 1} 
\end{array} \right] \, \left[ \begin{array}{c} r(t) \\ v(t) \\ \alpha(t) \end{array} \right] \, dt + 
\left[ \begin{array}{c} 0_{2 \times 3} \\ I_{3 \times 3} \end{array} \right] d\beta(t),                                                        
\end{align}
where the gravity related force term $G(x,t)$ and drag related force term $D(x,t)$ are given by
\begin{align}
G(x,t) & = -\frac{G m_0}{\|r(t)\|^3}, \\
D(x,t) & = -\beta_0 e^{\alpha} \exp\left\{ \frac{R_0 - \|r(t)\|}{H_0} \right\} \|v(t)\|.
\end{align}
The diffusion matrix for the Brownian motion $\beta(t)$ is 
\begin{equation}
Q(t) = \left[ \begin{array}{ccc} 2.4064 \cdot 10^{-5} & 0 & 0 \\
                                               0 & 2.4064 \cdot 10^{-5} & 0 \\
                                               0 & 0 & 1 \cdot 10^{-5} \end{array} \right].
\end{equation}
The effective diffusion matrix for this model is singular and the dynamic model can be written in the form of Equations (\ref{eq:singular_f1}) and (\ref{eq:singular_f2}). 
The values $\beta_0 = -0.59783$, $H_0 = 13.406$, $Gm_0 = 3.9860 \cdot 10^5$ and $R_0 = 6374$ are used as typical values for the parameters \cite{Julier2004} (see \cite{Leondes1981}). 

A radar located at $s = [s_x, s_y]^T$ periodically measures the range and bearing of the vehicle with 1 Hz frequency.
The measurement model is given by
\begin{align}
y_k & =\left[ \begin{array}{c} \| r(t_k) - s \| \\ \tan^{-1} \left( \frac{r_2(t_k) - s_y}{r_1(t_k) - s_x} \right) \end{array} \right] + v_k,
\end{align}
where the measurement noise $v_k$ is zero-mean Gaussian with covariance matrix
\begin{equation}
R_k = \left[ \begin{array}{cc} 1 \cdot 10^{-3} & 0 \\ 0 & 1.7 \cdot 10^{-3} \end{array} \right]. 
\end{equation}

The state trajectory and noisy measurements are simulated from $t_0 = 0\,$ to $t_K = 200\,$ using Euler--Maruyama discretization with time-step $\Delta t = 0.01$. 
The initial state is drawn from a Gaussian prior with mean and covariance given by
\begin{align}
m(t_0) & = \left[ \begin{array}{ccccc} 6500.4 & 349.14 & -1.8093 & -6.7967 & 0.6932 \end{array} \right]^T,  \\
P(t_0) & = \left[ \begin{array}{cc}  10^{-6}\cdot I_{4 \times 4} & 0_{4 \times 1} \\ 0_{1 \times 4} & 0 \end{array} \right]
 \label{eq:reentry_initial}.
\end{align}

For this model the computation of the Gaussian expectations needed in the GFGS and VGS is not possible analytically. 
The Gaussian expectations were computed using the EXT, CT, UT and G-H integration rules. 
For the VGS it was necessary to use the alternative formulation of CT2 and UT2 rules, since the CT and UT rules caused numerical problems and failure of the algorithm to converge. 
For G-H this was not a problem. 
The differential equations were solved using the standard 4-stage Runge--Kutta method with integration step of $0.1$. 
The initial mean and covariance for the smoothers are given by Equation (\ref{eq:reentry_initial}), where we used $m_5(t_0) = 0$ and $P_{5,5}(t_0) = 1$ for the unknown aerodynamic parameter. 

\begin{figure}[h!]
\centering
\includegraphics[scale=0.25]{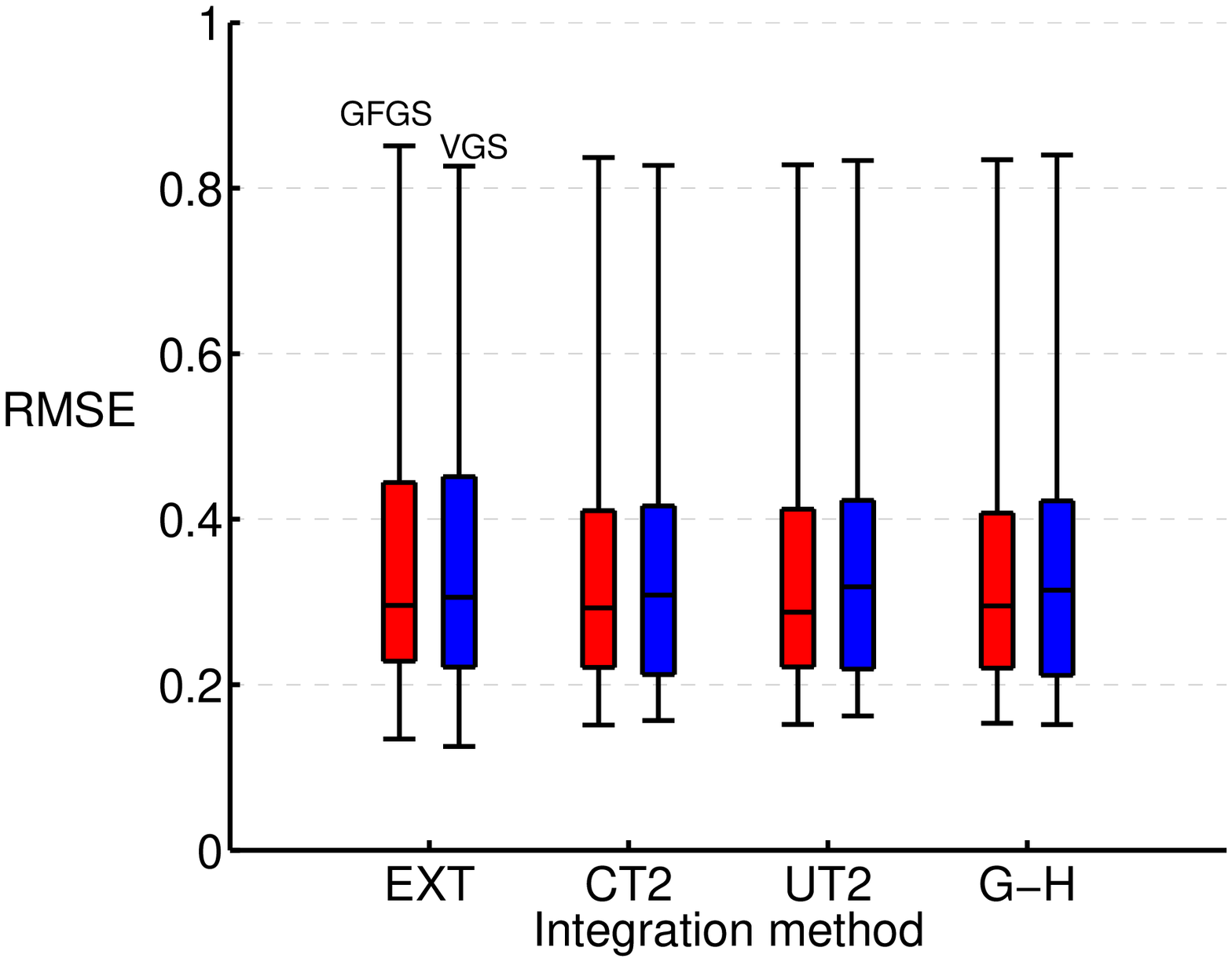} 
\includegraphics[scale=0.25]{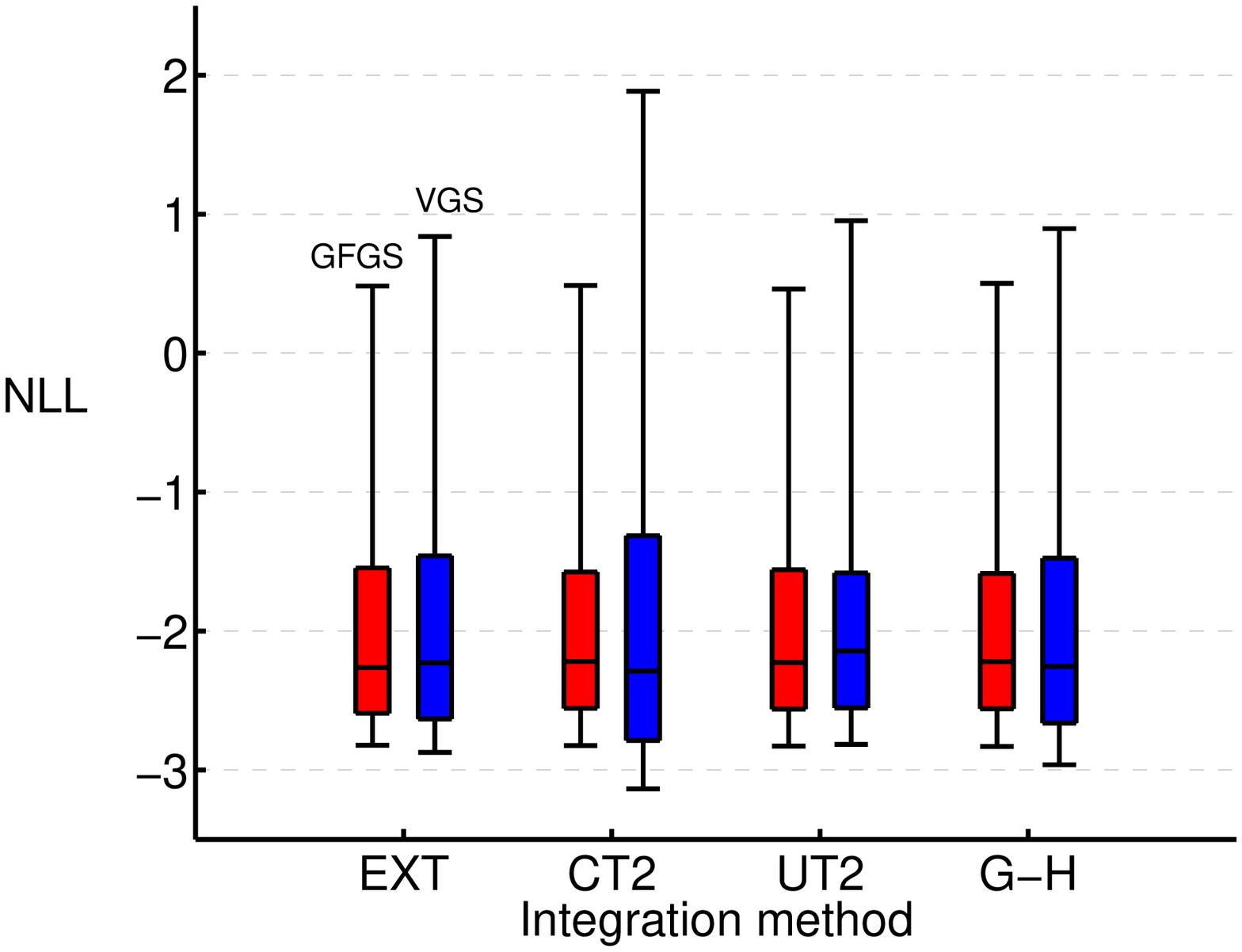} 
\includegraphics[scale=0.25]{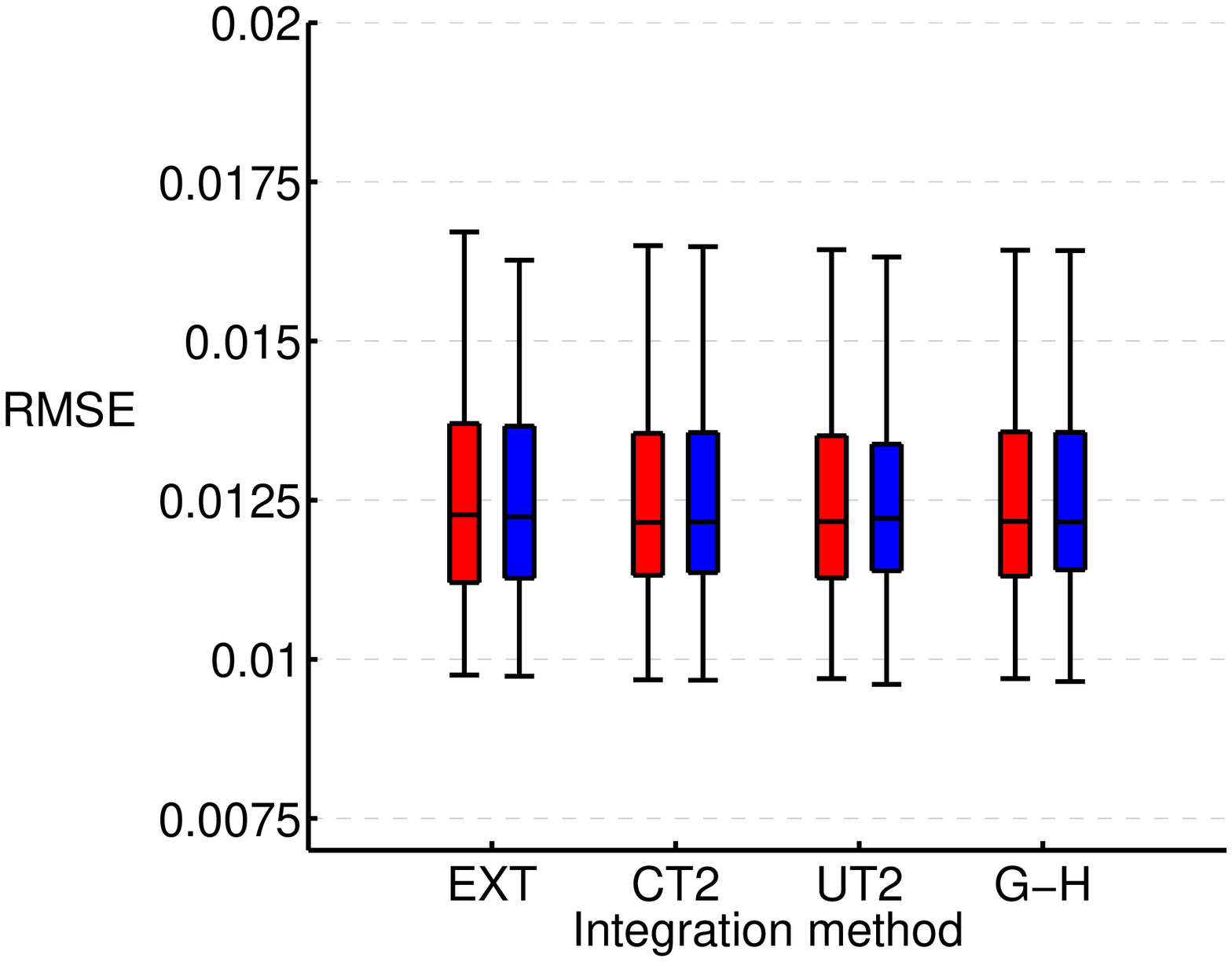}
\includegraphics[scale=0.25]{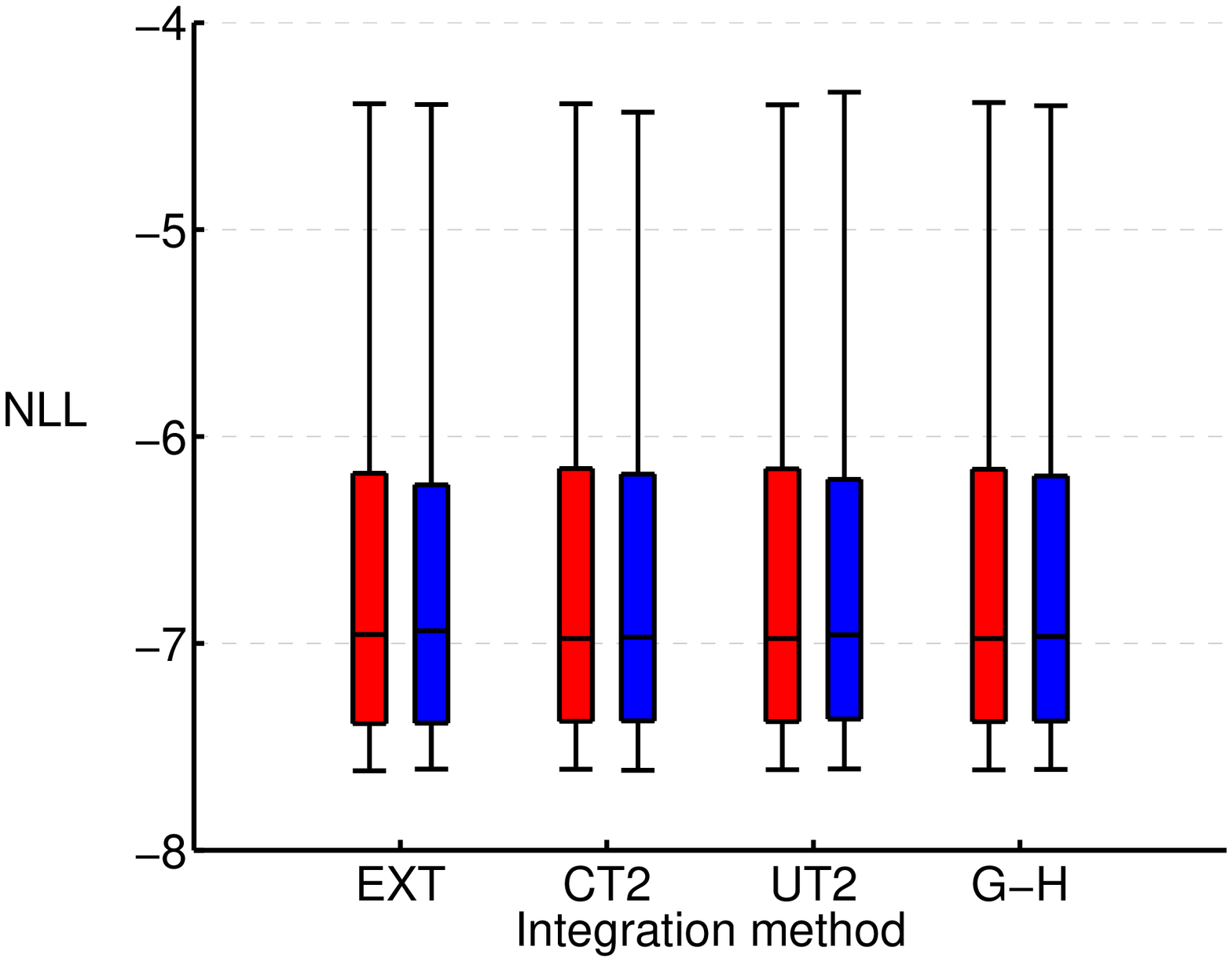}  
\includegraphics[scale=0.25]{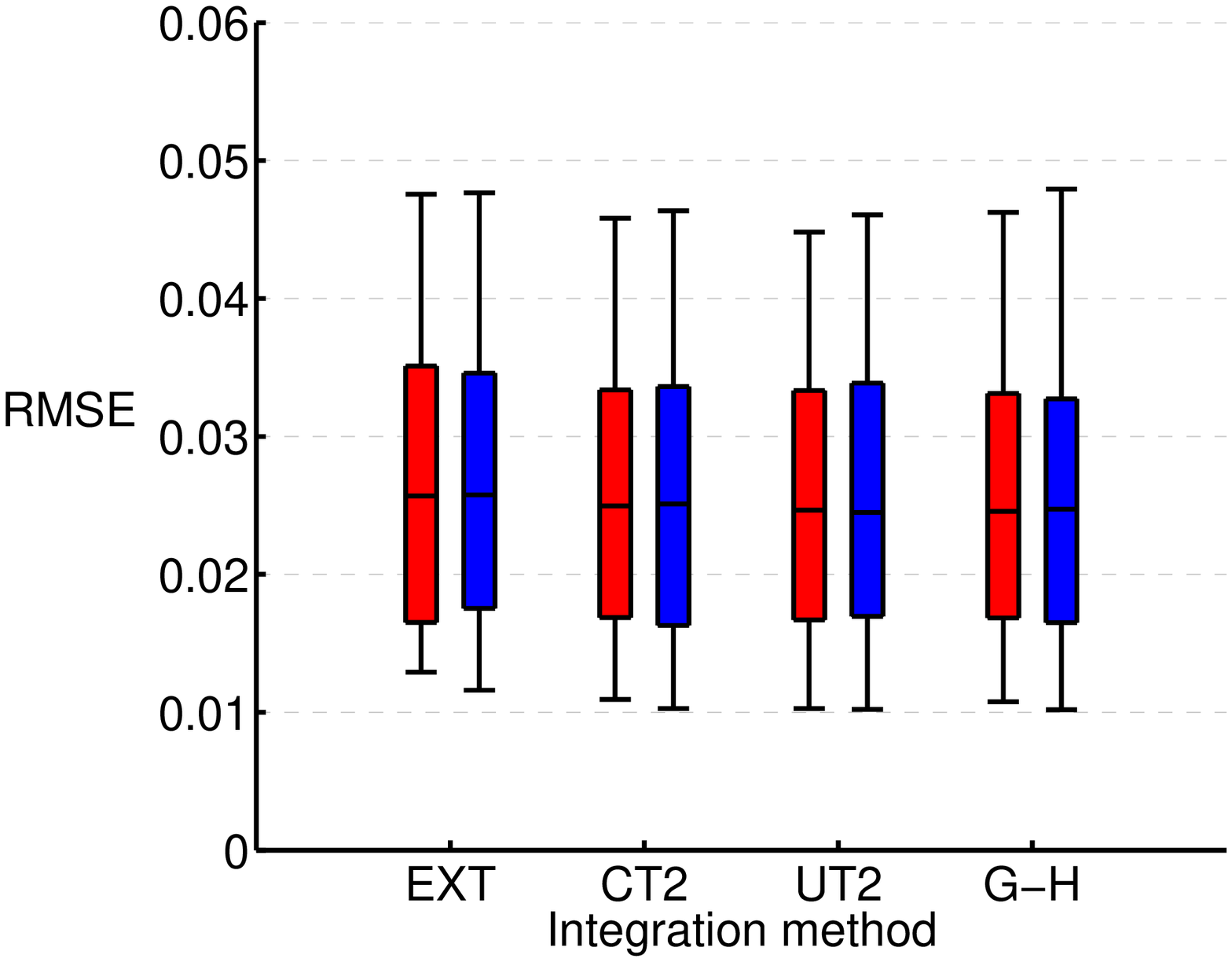}
\includegraphics[scale=0.25]{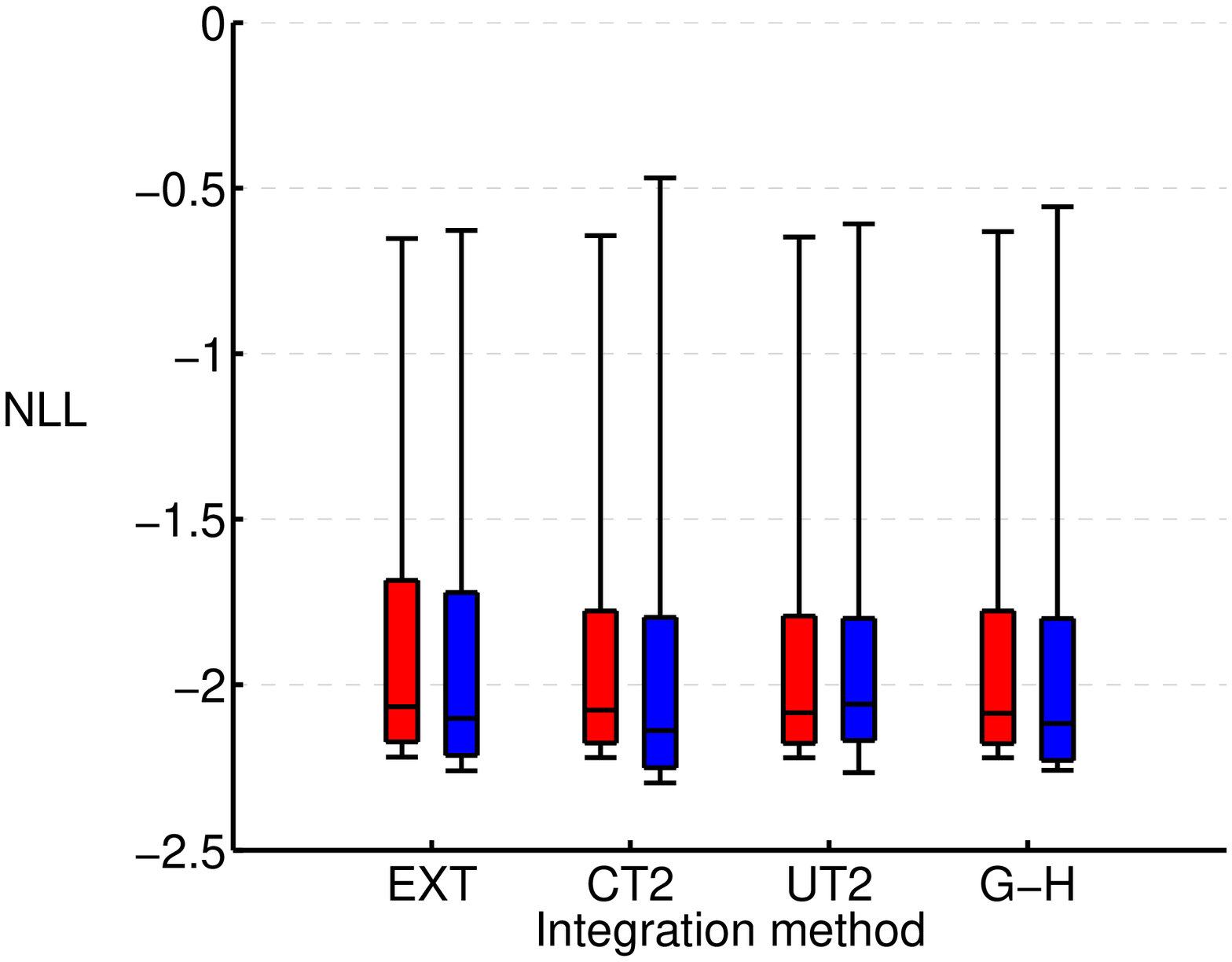}
\caption{The RMSE and NLL values for the Gaussian filtering based (red) and variational (blue) Gaussian smoothers using different methods to compute the Gaussian expectations. The top row shows the results for position, the middle row for velocity and the bottom row for the aerodynamic parameter. The boxplots show the 5\%, 25\%, 50\%, 75\% and 95\% quantiles.}
\label{fig:ex2_results} 
\end{figure}

The boxplots of RMSE and NLL results for 100 Monte Carlo simulations are shown in Figure \ref{fig:ex2_results}.
No clear difference can be seen between GFGS and VGS methods, or between the different methods for computing the Gaussian expectations.
The NLL values are slightly better for the GFGS, especially for the position. 
This is the result of the slightly underestimated variance for the position when using the VGS method. 
The more compact variance estimate of the VGS can also be seen from the mean 95\%-consistency results in Table \ref{table:ex2}.
On average, only 5 iterations were needed in the VGS before convergence. 
Also, only a small decrease of KL-divergence was observed, which explains the small difference between the methods. 

\begin{table}
\caption{The mean 95\%-consistencies for the Gaussian filtering based and variational Gaussian smoothers. }
\label{table:ex2}
\center
\begin{tabular}{r | c c c }
 & Position & Velocity & Parameter \\
 \hline
EXT-GFGS & 0.94 & 0.95 & 0.95  \\
EXT-VGS & 0.91 & 0.95 & 0.95 \\
CT-GFGS & 0.95 & 0.95 & 0.96  \\
CT2-VGS & 0.88 & 0.94 & 0.94 \\
UT-GFGS & 0.95 & 0.95 & 0.96  \\
UT2-VGS & 0.94 & 0.95 & 0.96 \\
GH-GFGS & 0.95 & 0.95 & 0.96  \\
GH-VGS & 0.92 & 0.95 & 0.95 
\end{tabular}
\end{table}

%
%

\section{Discussion and conclusions}

Compared to the Gaussian filtering based Gaussian smoother the variational Gaussian smoother is more complex to implement and is computationally heavier, since each iteration requires approximately the same amount of computations as one run of the filtering and smoothing equations. 
The Gaussian filtering based Gaussian smoother provides good initial conditions for the variational Gaussian smoother and could solve the problems with initialization mentioned in \cite{Vrettas2010}.
The variational Gaussian smoother provides Gaussian approximation that is optimal in the sense that it minimizes the Kullback--Leibler divergence of the approximating distribution with respect to the true distribution. 
However, the examples considered here show that this will not always improve the estimate with respect to other commonly used metrics. 

Using statistical linearization with respect to the filtering distribution in the variational Gaussian smoothing equations gives formally the Gaussian filtering based smoother as a special case. 
This suggests that for highly nonlinear systems, the variational Gaussian smoother could better capture the nonlinearities. 
This is demonstrated in the numerical experiment for the double well system, where the variational Gaussian smoother clearly improves the Gaussian filtering based smoother estimate for small measurement variances. 
However, no clear improvement was observed in the second numerical experiment for the reentry problem. 
The reason could be that the nonlinearities in the reentry system are not high enough to gain benefit from using the variational Gaussian smoother equations. 
A drawback of the variational Gaussian smoother is that it tends to underestimate the variance compared to the Gaussian filtering based smoother.  

For general nonlinear systems the Gaussian expectations can be computed using Taylor seires based linearization or standard sigma-point methods.
The Taylor series based variational Gaussian smoother linearizes the drift and measurement function with respect to the current mean estimate and is therefore similar in idea to the iterated Extended Kalman smoother \cite{Bell1994}.  
Using unscented transform and cubature rule based sigma-point methods in the variational Gaussian smoother resulted in some numerical problems. 
This is caused since the unscented and cubature methods are only accurate up to third order monomials, but even for linear systems some of the expectations are over fourth order polynomials. 
The numerical problems can be reduced by computing the Jacobians of the drift and measurement functions and using an alternative form for the expectations. 
Also, higher order unscented transform could be used for these expectations (see \cite{Julier2004, Lerner2002}).
Third order Gauss--Hermite integration rule worked well for both examples considered in this paper, but the computational cost is high especially for high dimensional systems. 
For some high dimensional systems the computational cost of the Gauss--Hermite method could be reduced by using Rao-Blackwellisation \cite{Closas2010}.

\section*{Acknowledgments}
The first author received financial support from the Tampere University of Technology Doctoral Programme in Engineering and Natural Sciences. 
The second author would like to thank the Academy of Finland.

\bibliographystyle{elsarticle-num}
\bibliography{bibliography_juha}

\appendix

\section{Computing the KL-divergence using Girsanov's theorem}
This section presents the derivation of the KL-divergence term in Equation (\ref{eq:KL_div_singular}) for the system with singular effective diffusion matrix. 
The KL-divergence term can be partitioned to \cite{ArchambeauOpper2011}
\begin{equation}
\text{KL}(\mathbb{Q}_X \, || \, \mathbb{P}_{X \cond Y}) = \text{KL}(\mathbb{Q}_X \, || \, \mathbb{P}_X) - \sum_{k=1}^K \E_q\left[ \ln p(y_k \, | \, x_k) \right], \label{eq:appA1}
\end{equation}
where $\mathbb{P}_X$ corresponds to the joint probability law of the stochastic processes $x_1(t)$ and $x_2(t)$ and $\mathbb{Q}_X$ to the joint probability law of the stochastic processes $s_1(t)$ and $s_2(t)$ defined by
\begin{align*}
\frac{dx_1}{dt} & = F_1(t) x, \\
 dx_2 & = f_2(x(t), t) + d \beta(t) \\
\frac{ds_1}{dt} & = F_1(t) s, \\
ds_2& = g_2(x(t), t) + d \beta(t),
\end{align*}
where $\beta(t)$ is Brownian motion with diffusion matrix $Q(t)$ with respect to measure $\mathbb{P}_X$.

The processes $s_1(t)$ and $s_2(t)$ are weak solutions to the original system under the measure $\mathbb{Q}_X$ that is defined through the Radon-Nikodym derivative \cite{SarkkaSottinen2007}
\begin{equation*}
\E \left[ \frac{d \mathbb{Q}_X}{d \mathbb{P}_X} \, | \, \mathrm{F}_t \right] = Z(t),
\end{equation*}
where $\mathrm{F}_t$ is the natural filtration of the Brownian motion $\beta(t)$ and
\begin{align*}
Z(t) & = \exp\left[ \int_0^t \left\{ f_2(s_1(t), s_2(t), t) - g_2(s_1(t), s_2(t), t) \right\}^T \, d \beta(t) \right. \nonumber \\
 & \left. -\frac{1}{2} \int_0^t \left\{ f_2(s_1(t), s_2(t), t) - g_2(s_1(t), s_2(t), t)\right\}^T Q^{-1}(t) \right. \nonumber 
 \\ & \left. \left\{ f_2(s_1(t), s_2(t), t) - g_2(s_1(t), s_2(t), t)^T\right\} \right]. 
\end{align*}
The KL-divergence term in the right side of Equation (\ref{eq:appA1}) is then given by
\begin{align*}
\text{KL}(\mathbb{Q}_X \, || \, \mathbb{P}_X) & = \E_{\mathbb{Q}_X} [- \ln Z(t)] \nonumber \\
 & = \frac{1}{2} \int_0^t \E_q \left[\left\{ f_2(s_1(t), s_2(t), t) - g_2(s_1(t), s_2(t), t)\right\}^T Q^{-1}(t)  \right. \nonumber 
 \\ & \left. \left\{ f_2(s_1(t), s_2(t), t) - g_2(s_1(t), s_2(t), t)\right\} \right].
\end{align*}
Inserting $g_2(s_1(t), s_2(t), t) = -A(t)s(t) + b(t)$ gives then the desired KL-divergence.

\section{Converting the Gaussian filtering based Gaussian smoother to the variational form}
Here we present the derivation of the variational form of the Gaussian filtering based Gaussian smoother. 
This is achieved by using the change of variables:
\begin{align*}
\lambda & = -P_f^{-1}(m_s-m_f), \quad \Psi  = -\frac{1}{2}\left(P_f^{-1}-P^{-1}_{s} \right).
\end{align*}
Computing the time derivatives of the new variables and inserting the filtering and smoothing differential equations (\ref{eq:dmf})-(\ref{eq:dPf}) and (\ref{eq:dms})-(\ref{eq:dPs}) gives
\begin{align*}
\frac{d}{dt}\lambda & = P_f^{-1}\left( \frac{d}{dt}P_f \right) P_f^{-1}(m_s-m_f) - P_f^{-1}\left( \frac{d}{dt}m^s - \frac{d}{dt}m_f \right) \\
& = P_f^{-1}\left(  \E_f[F_x(x)]P_f + P_f \E[F_x^T(x)] + Q \right)P_f^{-1}(m_s-m_f) \\
 &  - P_f^{-1}( \E_f[f(x)] + \E_f[F_x(x)](m_s-m_f) + QP_f^{-1}(m_s-m_f) - \E_f[f(x)]) \\
 & = -\E_f[F_x(x)]^T \lambda
\end{align*}
and 
\begin{align*}
\frac{d}{dt}\Psi & = \frac{1}{2}P_f^{-1}\left( \frac{d}{dt}P_f \right)P_f^{-1} - \frac{1}{2} P_s^{-1}\left(\frac{d}{dt}P_s \right)P_s^{-1} \\
 & = \frac{1}{2}P_f^{-1}( \E_f[F_x(x)]P_f + P_f \E[F_x(x)]^T + Q)P_f^{-1}  \\
  & -\frac{1}{2} P_s^{-1}(\E_f[F_x(x)]P_s + QP_f^{-1}P_s + P_s\E_f[F_x^T(x)] + P_sP_f^{-1}Q - Q)P_s^{-1} \\
  & = -\Psi\E_f[F_x(x)] - \E_f[F_x(x)]^T\Psi +2\Psi Q \Psi.
\end{align*}
Inserting $A(t) = -\E_f[F_x(x)] + 2 \Sigma(t) \Psi(t)$ to the above equations gives
\begin{align*}
\frac{d}{dt}\lambda(t) &= A^T(t) \lambda(t) - 2 \Psi(t) \Sigma(t) \lambda(t) \\
\frac{d}{dt}\Psi(t) &  = \Psi(t)A(t) +A^T(t)\Psi(t) - 2\Psi(t) \Sigma(t) \Psi(t).
\end{align*}

\section{Gradients with respect to mean and covariance}
In this section we derive the expressions for the gradients with respect to the mean vector $m$ and covariance matrix $P$ of a Gaussian expectation over a scalar function $e(x)$. 
These are used to form the sigma-point approximations for the Gaussian expectations in the variational Gaussian smoothing equations. 
Computing the gradients gives
\begin{align*}
\nabla_m \E[e(x)] & = \nabla_m \left[\int e(x) \text{N}(x \, | \, m, P) \, dx \right] = \int e(x) \nabla_m \text{N}(x, \, | \, m, P) \, dx \\ 
 & = \int e(x) \text{N}(x, \, | \, m, P)P^{-1}(x-m) \, dx  \\
 & = P^{-1}\E [e(x)(x-m)] = \E[\nabla_x e(x)]
\end{align*}
and
\begin{align*}
\nabla_P \E[e(x)] & = \nabla_P \left[ \int e(x) \text{N}(x \, | \, m, P) \, dx \right]  = \int e(x) \nabla_ P \text{N}(x \, | \, m, P) \, dx \\
 &  = \int e(x) \text{N}(x \, | \, m ,P)\frac{1}{2}\left[P^{-1}(x-m)(x-m)^T P^{-1} - P^{-1} \right] dx \\
 & = \frac{1}{2}P^{-1}\E[e(x)(x-m)(x-m)^T]P^{-1} - \frac{1}{2}\E[e(x)]P^{-1} \\
 & = \frac{1}{2}P^{-1}\E [\nabla_x e(x) (x-m)^T].
\end{align*}

\end{document}